\documentclass[11pt]{amsart}

\usepackage{amssymb}
\def\C{\mathbb{C}}
\def\N{\mathbb{N}}

\def\D{\mathbb D}

\def\K{\mathcal K}

\def\j{\mathcal J}

\def\s{\psi^{-1}}
\def\p{\varphi^{-1}}
\def\v{\varphi}
\def\l{\lambda}

\def\ov{\overline}
\def\lo{\longrightarrow}

\def\a{\alpha}
\def\b{\beta}
\def\g{\gamma}
\def\d{\displaystyle\sum}

\def\wi{\widetilde}
\def\e{equivalent}

\def\u{unitarily equivalent}

\def\ue{unitary equivalence}

 \newtheorem{thm}{Theorem}[section]
 \newtheorem{cor}[thm]{Corollary}
 \newtheorem{lem}[thm]{Lemma}
 \newtheorem{prop}[thm]{Proposition}
 
 \newtheorem{rem}[thm]{Remark}
 
\newcommand{\be}{\begin{equation}}
\newcommand{\ee}{\end{equation}}
\newcommand{\bea}{\begin{eqnarray}}
\newcommand{\eea}{\end{eqnarray}}
\newcommand{\Bea}{\begin{eqnarray*}}
\newcommand{\Eea}{\end{eqnarray*}}

\newcommand{\chd}{Cowen-Douglas }%
\newcommand{\dt}{$m$-tuple }%

\def\op{{operator~}}%
\def\ops{{operators~}}%
\def\rep{{representation~}}%
\def\HS{{Hilbert space}}%
\def\hs{{Hilbert space~}}%

\newcommand{\rk}{reproducing kernel }%
\newcommand{\Rk}{reproducing kernel}%

\newcommand{\bb}[1]{\mathbb #1}
\newcommand{\cl}[1]{\mathcal #1}
\newcommand{\alg}[1]{\mbox{${\mathcal A}(#1)$}}%
\def\F{\mbox{${\mathbb F}$}}
\newcommand{\kf}[1]{{K(\cdot,#1)}}%

\def\x{\mbox{{\boldmath $x$}}}%
\def\y{\mbox{{\boldmath $y$}}}%
\def\e{\mbox{{\boldmath $e$}}}%

\newcommand{\bm}[1]{{\mbox{\boldmath ${#1}$}}}

\newcommand{\inner}[2]{\langle #1,#2 \rangle }%

\newcommand{\qfor}{~\text{for}~}
\newcommand{\qand}{~\text{and}~}

\newcounter{cnt1}
\newcounter{cnt2}
\newcounter{cnt3}
\newcommand{\blr}{\begin{list}{$($\roman{cnt1}$)$}
 {\usecounter{cnt1} \setlength{\topsep}{0pt}
 \setlength{\itemsep}{0pt}}}
\newcommand{\bla}{\begin{list}{$($\alph{cnt2}$)$}
 {\usecounter{cnt2} \setlength{\topsep}{0pt}
 \setlength{\itemsep}{0pt}}}
\newcommand{\bln}{\begin{list}{$($\arabic{cnt3}$)$}
 {\usecounter{cnt3} \setlength{\topsep}{0pt}
 \setlength{\itemsep}{0pt}}}
\newcommand{\el}{\end{list}}

\begin{document}

\title[Irreducible homogeneous operators]{On the irreducibility of a class of \\homogeneous operators}
\author[Misra]{\!Gadadhar Misra} 
\author[Shyam Roy]{Subrata Shyam Roy}
\address{ Indian Statistical Institute \\
R. V. College Post Office\\ Bangalore -- 560 059\\ India} 
\email[Gadadhar Misra]{gm@isibang.ac.in} 
\email[Subrata Shyam Roy]{ssroy@isibang.ac.in}




\sloppy
\date{}


\thanks{The research of the first author was supported in part by a
grant from the DST - NSF  Science and Technology Cooperation Programme.  
The second author was supported by the Indian Statistical Institute.}


\begin{abstract}
In this paper we construct a class of homogeneous Hilbert modules over
the disc algebra $\mathcal{A}(\D)$ as quotients of certain natural
modules over the function algebra $\mathcal{A}(\D^2)$.  
These quotient modules are described using the jet construction for Hilbert
modules. We show that the quotient modules obtained this way, belong to the 
class ${\mathrm B}_k(\D)$ and that they are mutually inequivalent, irreducible
and homogeneous.  
\end{abstract}
\maketitle
\section{Introduction}
Let $\cl{M}$ be a Hilbert space.
All Hilbert spaces in this paper will be assumed to be
complex and separable.
Let $\cl{A}(\Omega)$ be the natural function algebra
consisting of functions holomorphic in a neighborhood of
the closure $\bar{\Omega}$ of some open, connected and bounded subset
$\Omega$ of $\C^m$. The
Hilbert space $\cl{M}$ is said to be a {\em Hilbert module\/} over
$\cl{A}(\Omega)$ if $\cl{M}$ is a module over $\cl{A}(\Omega)$ and
$$
\|f \cdot h\|_{\cl M} \le C \|f\|_{\cl A(\Omega)} \|h\|_{\cl M}
\qfor f \in \cl A(\Omega) \qand h \in \cl M,
$$
for some positive constant $C$ independent of $f$ and $h$.
It is said to be {\em contractive \/} if we also have
$C \leq 1$.

Fix an inner product on the algebraic tensor product $\alg{\Omega}
\otimes \C^n$. Let the completion of $\alg{\Omega} \otimes \C^n$
with respect to this inner product be the Hilbert space $\mathcal
M$. Assume that the module action
$$\alg{\Omega} \times \alg{\Omega} \otimes \C^n\to \alg{\Omega} \otimes \C^n$$ extends
continuously to $\alg{\Omega} \times \cl{M} \to \cl{M}$. With very
little additional assumption on $\mathcal{M}$, we obtain a {\em
quasi-free} Hilbert module (cf. \cite{quasifree}).

The simplest family of modules over $\mathcal{A}(\Omega)$ corresponds
to evaluation at a point in the closure of $\Omega$.  For $\underline{z}$ in
the closure of $\Omega$, we make the one-dimensional Hilbert space
$\C$ into the Hilbert module $\C_{\underline{z}}$, by setting $\varphi v =
\varphi(z) v$ for $\varphi$ $\in \mathcal{A}(\Omega)$ and $v \in \C$.
Classical examples of contractive Hilbert modules are the Hardy and Bergman
modules over the algebra $\mathcal{A}(\Omega)$.

Let $G$ be a locally compact second countable group acting
transitively on $\Omega$.  Let us say that the module $\cl{M}$ over
the algebra $\mathcal A(\Omega)$ is homogeneous if $\varrho(f \circ
\varphi)$ is unitarily equivalent to $\varrho(f)$ for all $\varphi \in
G$.  Here $\varrho: \mathcal A(\Omega) \to \mathcal B(\cl{M})$ is the
homomorphism of the algebra $\mathcal A(\Omega)$ defined by
$\varrho(f) h := f \cdot h$ for $f\in \alg{\Omega}$ and $h\in \cl{M}$.
It was shown in \cite{MSJOT90} that if the module $\cl{M}$ is irreducible and homogeneous then
there exists a {\em projective unitary representation} $U: G \to \mathcal U(\cl{M})$ such that
$$
U_\varphi^* \varrho(f) U_\varphi = \varrho(f\cdot \varphi),\; f
\in \alg{\Omega},\; \varphi\in G,
$$ where $(f\cdot \varphi)(w) = f(\varphi\cdot w)$ for $w\in \Omega$.

A $*$ - homomorphism $\varrho$ of a $C^*$ - algebra $\mathcal C$
and a unitary group representation $U$ of $G$ on the Hilbert space
$\mathcal M$ satisfying the condition as above were first studied
by Mackey and were called {\em Systems of Imprimitivity}.  Mackey
proved the Imprimitivity theorem which sets up a correspondence
between induced representations of the group $G$ and the Systems
of Imprimitivity.  The notion of homogeneity is obtained by
compressing the systems of imprimitivities, in the sense of
Mackey, to a subspace $\mathcal N$ of $\mathcal M$ and then
restricting to a subalgebra of the $C^*$ - algebra $\mathcal C$
(cf.\cite{impr}).  However, it is not clear if the notion of
homogeneity is in some correspondence with holomorphically induced
representations, at least when the module $\cl{M}$ is assumed to
be in $\mathrm B_k(\Omega)$.

An alternative description, in the particular case of the disc may be
useful. The group of bi-holomorphic automorphisms M\"{o}b of the unit disc is
$\{\varphi_{\theta, \alpha}:\, \theta\in [0, 2\pi) \qand \alpha
\in \D\}$, where
\begin{equation} \label{mobius map} \varphi_{\theta,\alpha}(z) = \e^{i \theta}\frac{z-
\a}{1-\bar{\alpha}z},~~z\in \D.
\end{equation}
As a topological group (with the
topology of locally uniform convergence) it is isomorphic to PSU($1,1$) and to
PSL($2,\mathbb{R}$).

An operator $T$ from a Hilbert space into itself is said to
be {\em homogeneous} if $\varphi(T)$ is unitarily equivalent to $T$
for all $\varphi$ in M\"{o}b which are analytic on the spectrum of
$T$. The spectrum of a homogeneous operator $T$ is either the unit
circle $\mathbb{T}$ or the closed unit disc $\bar{\mathbb{D}}$, so
that, actually, $\varphi(T)$ is unitarily equivalent to $T$ for all
$\varphi$ in M\"{o}b. We say that a projective unitary representation
$U$ of M\"{o}b is {\em associated} with an operator $T$ if
\[
\varphi(T)=U_\varphi^{\ast}TU_\varphi
\]
for all $\varphi$ in M\"{o}b.  We have already pointed out that
if $T$ is irreducible then it has an associated representation
$U$.  It is not hard to see that $U$ is uniquely determined upto
unitary equivalence.

Many examples (unitarily inequivalent) of homogeneous operators
are known \cite{BMIAS01}. Since the direct sum (more generally
direct integral) of two homogeneous operators is again
homogeneous, a natural problem is the classification (up to
unitary equivalence) of {\em atomic homogeneous operators}, that
is, those homogeneous operators which can not be written as the
direct sum of two homogeneous operators. In this generality, this
problem remains unsolved. However, the irreducible homogeneous
operators in the Cowen-Douglas class $\mathrm B_1(\D)$ and
$\mathrm B_2(\D)$ have been classified (cf. \cite{hom} and
\cite{wil}) and all the scalar shifts (not only the irreducible
ones) which are homogeneous are known \cite[List 4.1, page
312]{shift}. Some recent results on classification of 
homogeneous bundles are in \cite{psl2r} and \cite{KMCompt}.

Clearly, irreducible homogeneous operators are atomic.
Therefore, it is important to understand when a homogeneous
operator is irreducible.

There are only two examples of atomic homogeneous
operators known which are not irreducible.  these
are the multiplication operators -- by the respective co-ordinate functions --
on the Hilbert spaces $L^{2}(\mathbb{T)}$ and $L^{2}(\mathbb{D)}$.
Both of these examples happen to be normal operators.
We do not know if all atomic homogeneous operators possess an associated
projective unitary representation. However, to every homogeneous operator
in ${\mathrm B}_k(\D)$, there exist an associated representation of the
universal covering group of M\"{o}b \cite[Theorem 4]{KMCompt}.

It turns out an irreducible homogeneous operator in $\mathrm
B_2(\D)$ is the compression of the tensor product of two
homogeneous operators from $\mathrm B_1(\D)$ (cf. \cite{BMIAS01})
to a suitable invariant subspace. In the language of Hilbert
modules, this is the statement that every homogeneous module in
$\mathrm B_2(\D)$ is obtained as quotient of the tensor product of
two homogeneous modules in $\mathrm B_1(\D)$ by the sub-module of
functions vanishing to order $2$ on $\triangle \subseteq \D^2$.
However, beyond the case of rank $2$, the situation is more
complicated. The question of classifying homogeneous operators in
the class $\mathrm B_k(\D)$ amounts to classifying holomorphic and
Hermitian vector bundles of rank $k$ on the unit disc which are
homogeneous. Classification problems such as this one are well
known in the representation theory of locally compact second
countable groups. However, in that context, there is no Hermitian
structure present which makes the classification problem entirely
algebraic.  A complete classification of homogeneous operators in
$\mathrm B_k(\D)$ may still be possible using techniques from the
theory of unitary representations of the M\"{o}bius group. Leaving
aside, the classification problem of the homogeneous operators in
${\mathrm B}_k(\D)$, we show that the  ``generalized Wilkins
examples" (cf. \cite{BMIAS01}) are irreducible.

If one considers a bounded symmetric domain in $\C^m$, the
classification question probably is even more complicated
(cf. \cite{BMJFA96}, \cite{Arazy-Z}). Here part of the difficulty lies
in the fact that no classification of the irreducible unitary
representations of the group ${\rm Aut}(\Omega)$, the bi-holomorphic
automorphism group of $\Omega$, is known.

In the following section, we discuss reproducing kernels for a functional Hilbert
space on a domain $\Omega \subseteq \C^m$ and the \dt of multiplication operators
$\bm{M}$ of multiplication by coordinate functions.  Although, our applications to the question
of irreducibility is only for the multiplication operator $\bm M$ on a functional
Hilbert space based on the unit disc $\D$, the more general discussion of this section 
is not any simpler in the one variable case.  

In  subsection \ref{CD}, we explain the realization of a \dt of operators $\bm T$ in the  
class $\mathrm B_k(\Omega)$ as the adjoint of a \dt of multiplication operators $\bm M$ on 
a Hilbert space of holomorphic functions, on the bounded connected open set 
$\Omega^*:=\{w \in \mathbb C^m: \bar{w} \in \Omega\}$,  possessing a reproducing kernel $K$.  
We point out, as in \cite{C-S},  that  the normalized kernel $\wi{\wi K}$ obtained from the 
kernel $K$  by requiring that  $\wi{\wi K}(z,0)=1$, for all $z\in \Omega$, determines the uniatry 
equivalence class of the \dt $\bm T$.  We then obtain a criterion for the irreducibility of the \dt 
$\bm T$ in terms of the normalized kernel $\wi{\wi K}$.  Roughly speaking, this says that the 
\dt of operators is irreducible if and only if the coefficients, in the pwer series expansion of   $\wi{\wi K}$,
are simultaneously irreducible.   Following, \cite{DMV} and \cite{DMTrans}, we describe the 
jet construction for Hilbert modules and discuss some examples.  

In section \ref{secmult}, we show that if $\mathcal H$ is a Hilbert space of holomorphic functions, on a 
bounded connected open set $\Omega$ and possesses a reproducing kernel $K$ then it admits a natural 
multiplier representation of the automorphism group  of $\Omega$ if $K$ is 
{\em quasi-invariant}.  We show that if $K$ is quasi-inavariant, then the corresponding multiplier 
representation intertwines $\bm M$ and $\varphi(\bm M)$, that is, the \dt of multiplication operators 
$\bm M$ is {\em homogeneous}.   

Our main results on irreducibility of certain class of homogeneous operators is in Section 4. 
The kernel $B^{(\alpha, \beta)}(z,w) = (1-z_1\bar{w}_1)^{-\alpha} (1-z_2\bar{w}_2)^{-\beta}$, $z=(z_1,z_2),\, w=(w_1,w_2) \in \D^2$, determines a Hilbert module over the function algebra $\mathcal A(\D^2)$. 
We recall the computation of a  matrix valued  kernel on the unit disc $\D$ using the 
jet construction for this Hilbert module which consists of holomorphic functions on the 
unit disc $\D$ taking values in $\C^n$.  The multiplication operator on this Hilbert space is then shown to be irreducible by 
checking that all of the coefficients of the ``normalized''  matrix valued kernel, obtained from the jet construction, cannot be simultaneously reducible.  

In section 5, we show that the kernel  obtained from the jet construction is quais-invariant and consequently, the 
corresponding multiplication operator is homogeneous.  This proof involves the verification of a cocycle identity, 
which in turn, depends on a beutiful identity involving binomial coefficients. 

Finally, in section 6, we discuss some examples arising from the jet construction applied to a certain natural family of 
Hilbert modules over the algebra $\mathcal A(\D^3)$.  A more systematic study of such examples 
is to be found in \cite{MSubCurv}. 

\section{Reproducing Kernels and the Cowen-Douglas class}
\subsection{\sf Reproducing kernel} \label{rkf}
Let ${\mathcal L}(\F)$ be the Banach space of all linear transformations on
a Hilbert space $\F$ of dimension $n$ for some $n \in \N$. Let
$\Omega\subset \bb{C}^m$ be a bounded open connected set.  A function
$K:\Omega\times\Omega \to {\mathcal{L}}(\F)$, satisfying
\begin{equation} \label{existence reprod}
\sum_{i,j=1}^p \inner{K(w^{(i)},w^{(j)})\zeta_j}{\zeta_i}_{\mathbb F}
~\geq~0 ,~~w^{(1)},\ldots,w^{(p)}\in \Omega, ~~\zeta_1,\ldots,\zeta_p
\in \F,~ p > 0
\end{equation}
is said to be a {\em non negative definite (nnd) kernel} on $\Omega$.
Given such an nnd kernel $K$ on $\Omega$, it is easy to construct a \hs
$\mathcal{H}$ of functions on $\Omega$ taking values in $\F$ with the
property
\begin{equation} \label{reproducing property}
\inner{f(w)}{\zeta}_{\mathbb F} = \inner{f}{\kf{w}\zeta}, \qfor w\in
\Omega,~\zeta\in \F, \qand f\in \mathcal{H}.
\end{equation}
The \hs $\mathcal{H}$ is simply the completion of the linear span of all
vectors of the form $\cl{S} = \{\kf{w}\zeta$, $w\in \Omega$, $\zeta\in
\F\}$, where the inner product between two of the vectors from $\cl{S}$ is
defined by
\begin{equation}\label{nndinner}
\inner{\kf{w}\zeta}{\kf{w^\prime}\eta} =
\inner{K(w^\prime,w)\zeta}{\eta}, \qfor \zeta,\eta\in\F,
\qand w,w^\prime \in \Omega,
\end{equation}
which is then extended to the linear span $\cl{H}^{\circ}$ of the set $\cl{S}$.
This ensures the reproducing property \mbox{(\ref{reproducing property})} of
$K$ on $\cl{H}^{\circ}$.
\begin{rem} \label{nndkernel}
We point out that although the kernel $K$ is required to be merely
{\em nnd}, the equation \mbox{(\ref{nndinner})} defines a positive
definite sesqui-linear form. To see this, simply note that
$|\inner{f(w)}{\zeta}| = |\inner{f}{\kf{w}\zeta}|$ which
is at most $\|f\|\inner{K(w,w)\zeta}{\zeta}^{1/2}$ by
the Cauchy - Schwarz inequality. It follows that if $\|f\|^2 = 0$
then $f=0$.
\end{rem}
Conversely, let $\mathcal{H}$ be any Hilbert space of functions on
$\Omega$ taking values in $\F$.  Let $e_w : \mathcal{H} \to
\F$ be the evaluation functional defined by $e_w(f) =
f(w)$, $w\in \Omega$, $f\in \mathcal{H}$.  If $e_w$
is bounded for each $w \in \Omega$ then it admits a bounded
adjoint $e_w^*: \F \to \cl{H}$ such that $\inner{e_w
f}{\zeta} = \inner{f}{e_w^* \zeta}$ for all $f\in \cl{H}$ and
$\zeta\in \F$.  A function $f$ in $\cl{H}$ is then orthogonal to
$e_w^*(\cl{H})$ if and only if $f=0$.  Thus $f = \sum_{i=1}^p
e_{w^{(i)}}^*(\zeta_i)$ with $w^{(1)},
\ldots,w^{(p)}\in \Omega,~\zeta_1,\ldots,\zeta_p \in \F,\qand
p > 0,$ form a dense set in $\cl{H}$.   Therefore we have
$$
\|f\|^2 = \sum_{i,j=1}^p \inner{e_{w^{(i)}}e_{w^{(j)}}^*\zeta_j}
{\zeta_i},
$$
where $f=\sum_{i=1}^n e_{w^{(i)}}^*(\zeta_i), ~w^{(i)}
\in \Omega, ~\zeta_i\in \cl{F}$. Since $\|f\|^2 \geq 0$, it
follows that the kernel $K(z,w) = e_ze_w^*$ is
non-negative definite as in \mbox{(\ref{existence reprod})}.  It
is clear that $K(z,w)\zeta \in \mathcal{H}$ for each
$w\in \Omega$ and $\zeta\in \F$, and that it has the
reproducing property \mbox{(\ref{reproducing property})}.
\begin{rem} \label{nonsingker}
If we assume that the evaluation functional $e_w$ is surjective then the
adjoint $e_w^*$ is injective and it follows that
$\inner{K(w,w)\zeta}{\zeta} > 0$ for all non-zero vectors $\zeta\in \F$.
\end{rem}

There is a useful alternative description of the reproducing
kernel $K$ in terms of the orthonormal basis $\{e_k: k \geq 0\}$
of the Hilbert space $\cl{H}$. We think of the vector $e_k(w)
\in \F$ as a column vector for a fixed $w\in \Omega$ and let
$e_k(w)^*$ be the row vector $(\overline{e_k^1(w)}, \ldots ,
\overline{e_k^n(w)})$. We see that
\begin{eqnarray*}
\inner{K(z,w)\zeta}{\eta} &=&
\inner{\kf{w}\zeta}{\kf{z}\eta}\nonumber\\
&=& \sum_{k=0}^\infty \inner{\kf{w}\zeta}{e_k} \inner{e_k}{\kf{z}\eta}
\nonumber\\
&=& \sum_{k=0}^\infty\overline{\inner{{e_k(w)}}{\zeta}}\inner{e_k(z)}{\eta}
\nonumber\\
&=&\sum_{k=0}^\infty \inner{e_k(z) e_k(w)^* \zeta}{\eta},
\end{eqnarray*}
for any pair of vectors $\zeta, \eta \in \F$.  Therefore, we have the following very
useful representation for the reproducing kernel $K$:
\begin{equation} \label{sumonb}
K(z,w) = \sum_{k=0}^\infty e_k(z)  e_k(w)^*,
\end{equation}
where $\{e_k: k \geq 0\}$ is any orthonormal basis in $\cl{H}$.

\subsection{\sf The Cowen-Douglas class} \label{CD}
Let $\bm T = (T_1, \ldots , T_m)$ be a d-tuple of commuting bounded linear \ops
on a separable complex \hs $\cl{H}$. Define the \op $D_{\mathbf T}:{\mathcal
H}\to {\mathcal H}\oplus  \cdots \oplus {\mathcal H}$ by $D_{\mathbf T}(x) = (T_1 x,
\ldots, T_m x)$, $x\in {\mathcal H}$.  Let $\Omega$ be a bounded domain in
$\C^m$.   For $w = (w_1, \ldots , w_m) \in \Omega$, let $\bm T-w$ denote the \op tuple
$(T_1-w_1, \ldots ,T_m - w_m)$.  Let $n$ be a positive integer.
\begin{def} \label{C-D def}
The \dt $\bm T$ is said to be in the \chd class ${\mathrm B}_n(\Omega)$ if
\begin{enumerate}
\item $\mbox{\rm ran}~D_{\mathbf T - w}$ is closed for all $w\in \Omega$
\item $\mbox{\rm span}~\{ \ker D_{\mathbf T - w}: w \in \Omega\}$ is
dense in ${\mathcal H}$
\item $\dim \ker D_{\mathbf T - w}= n$  for all $w \in \Omega$.
\end{enumerate}
\end{def}
This class was introduced in \cite{Bolyai}.  The case of a single \op was
investigated earlier in the paper \cite{C-D}.  In this paper, it is
pointed out that an operator $T$ in ${\mathrm B}_1(\Omega)$ is unitarily
equivalent to the adjoint of the multiplication \op $M$ on a \rk \HS,
where $(Mf)(z)=zf(z)$.  It is not very hard to see that, more generally, 
a \dt $\bm T$ in ${\mathrm B}_n(\Omega)$ is unitarily equivalent to the 
adjoint of the \dt of multiplication \ops
$\bm M = (M_1, \ldots ,M_m)$ on a reproducing kernel \hs \cite{C-D} and
\cite[Remark 2.6 a) and b)]{C-S}.  Also, Curto and Salinas \cite{C-S}
show that if certain conditions are imposed on the \rk then the
corresponding adjoint of the \dt of multiplication \ops belongs to the class
${\mathrm B}_n(\Omega)$.

To a \dt $\bm T$ in ${\mathrm B}_n(\Omega)$, on the one hand, one may
associate a holomorphic Hermitian vector bundle $E_{\mathbf T}$ on 
$\Omega$ (cf. \cite{C-D}), while on the other hand, one may associate 
a normalized \rk $K$ (cf. \cite{C-S}) on a suitable sub-domain of 
$\Omega^*=\{w \in \C^m: \bar{w}\in \Omega\}$.  It is possible to 
answer a number of questions regarding the \dt of \ops $\bm T$ using 
either the vector bundle or the \Rk.  For instance, in the two papers 
\cite{C-D} and \cite{Bolyai}, Cowen and Douglas show that the curvature
of the bundle $E_{\mathbf T}$ along with a certain number of derivatives forms
a complete set of unitary invariants for the operator $\bm T$ while
Curto and Salinas \cite{C-S} establish that the unitary equivalence
class of the normalized kernel $K$ is a complete unitary invariant for 
the corresponding \dt of multiplication operators.  Also, in \cite{C-D}, 
it is shown that a single operator in ${\mathrm B}_n(\Omega)$ is reducible 
if and only if the associated holomorphic Hermitian vector bundle admits 
an orthogonal direct sum decomposition.

We recall the correspondence between a \dt of operators in the class
${\mathrm B}_n(\Omega)$ and the corresponding \dt of multiplication
operators on a \rk Hilbert space on $\Omega$.

Let $\bm T$ be a \dt of operators in ${\mathrm B}_n(\Omega)$. Pick $n$ linearly
independent vectors $\gamma_1(w), \ldots,\gamma_n(w)$
in $\ker D_{\mathbf T-w}$, $w\in \Omega$. Define a map
$\Gamma: \Omega \to {\mathcal L}(\F,{\mathcal H})$ by
$\Gamma(w)\zeta = \sum_{i=0}^n \zeta_i\gamma_i(w)$, where $\zeta=
(\zeta_1, \ldots, \zeta_n)\in \F$, $\dim\F = n$.  It is shown in
\cite[Proposition 1.11]{C-D} and \cite[Theorem 2.2]{C-S} that it is possible to
choose $\gamma_1(w), \ldots ,\gamma_n(w)$, $w$ in some
domain $\Omega_0 \subseteq \Omega$, such that
$\Gamma$ is holomorphic on $\Omega_0$.
Let ${\mathcal A}(\Omega, \F)$ denote the
linear space of all $\F$ - valued holomorphic functions on $\Omega$.
Define $U_\Gamma: {\mathcal H}\to {\mathcal A}(\Omega_0^*, \F)$ by
\begin{equation} \label{general construction}
(U_\Gamma x)(w) = \Gamma(w)^* x, ~~x\in {\mathcal H},
~w\in \Omega_0.
\end{equation}
Define a sesqui-linear form on ${\mathcal H}_\Gamma =
\mbox{ran}~U_\Gamma$ by $\inner{U_\Gamma f}{U_\Gamma g}_\Gamma =
\inner{f}{g}$, $f,g\in{\mathcal H}$. The map $U_\Gamma$ is linear
and injective.  Hence ${\mathcal H}_\Gamma$ is a \hs of
$\F$-valued holomorphic functions on $\Omega_0^*$ with inner
product $\inner{\cdot}{\cdot}_\Gamma$ and $U_\Gamma$ is unitary.
Then it is easy to verify the following (cf. \cite[Remarks
2.6]{C-S}).
\begin{enumerate}
\item[a)] $K(z,w) = \Gamma(\bar{z})^* \Gamma (\bar{w})$,
$z,w \in \Omega_0^*$ is the \rk for the \hs ${\mathcal
H}_\Gamma$. \item[b)] $M_i^* U_\Gamma = U_\Gamma T_i$, where $(M_i
f)(z) = z_i f(z)$, $z=(z_1, \ldots , z_m) \in \Omega$.
\end{enumerate}
\begin{def} \label{normalised}
An nnd kernel $K$ for which $K(z,w_0) = I$ for all
$z\in \Omega_0^*$ and
some $w_0 \in \Omega$ is said to be {normalized at $w_0$}.
\end{def}

For $1\leq i\leq m$, suppose that the \ops $M_i:{\mathcal H} \to
{\mathcal H}$ are bounded.  Then it is easy to verify that for each fixed
$w\in \Omega$, and $1\leq i \leq m$,
\begin{equation} \label{eigenspace}
M_i^* K(\cdot, w)\eta = \bar{w}_i K(\cdot, w)\eta
\qfor \eta \in \F.
\end{equation}
Differentiating \mbox{ (\ref{reproducing property})}, we also obtain the following
extension of the reproducing property:
\begin{equation} \label{ext reprod}
\inner{(\partial_i^jf)(w)}{\eta}=\inner{f}{\bar{\partial}^j_i
K(\cdot,w)
\eta}~~\mbox{for~}1\leq i\leq m,~~j\geq 0,~w\in \Omega,~\eta\in \F,~f\in
{\mathcal H}.
\end{equation}
Let $\bm{M} = (M_1, \ldots , M_m)$ be the commuting $m$ - tuple of multiplication
operators and let $\bm{M}^*$ be the $m$ - tuple $(M_1^*,\ldots ,M_m^*)$. It then
follows from \mbox{(\ref{eigenspace})} that the eigenspace of the $m$ - tuple
$\bm{M}^*$ at $w \in \Omega^*\subseteq \C^m$ contains the $n$-dimensional subspace
${\mathrm ran}\,K(\cdot,\bar{w}) \subseteq \cl{H}$.

One may impose additional conditions on $K$ to ensure that $\bm{M}$ is in ${\mathrm
B}_n(\Omega^*)$. Assume that $K(w,w)$ is invertible for $w\in \Omega$. Fix $w_0 \in
\Omega$ and note that $K(z, w_0)$ is invertible for $z$ in some neighborhood
$\Omega_0 \subseteq \Omega$ of $w_0$. Let $K_{\rm res}$ be the restriction of $K$ to
$\Omega_0\times\Omega_0$. Define a kernel function $K_0$ on $\Omega_0$ by
\begin{equation} \label{normalised kernel}
K_0(z,w) = \varphi(z) K(z,w)\varphi(w)^*,~ z,w\in \Omega_0,
\end{equation}
where $\varphi(z) =  K_{\rm res}(w_0,w_0)^{1/2}K_{\rm res}(z, w_0)^{-1}$. The kernel
$K_0$ is said to be {\em normalized} at $0$ and is characterized by the property
$K_0(z,w_0)=I$ for all $z\in \Omega_0$. Let  $\bm{M}_0$ denote the \dt of
multiplication operators on the Hilbert space $\cl H$. It is not hard to establish
the unitary equivalence of the two $m$ - tuples $\bm{M}$ and $\bm{M}_0$ as in (cf.
\cite[Lemma 3.9 and Remark 3.8]{C-S}). First, the restriction map $res: f \to
f_{res}$, which restricts a function in ${\mathcal H}$ to $\Omega_0$ is a unitary
map intertwining the \dt $\bm{M}$ on ${\mathcal H}$ with the \dt $\bm{M}$ on
${\mathcal H}_{\rm res}= {\rm ran}~res$. The \hs ${\mathcal H}_{\rm res}$ is a \rk
\hs with \rk $K_{\rm res}$. Second, suppose that the $m$ - tuples $\bm{M}$ defined on
two different \rk Hilbert spaces ${\mathcal H}_1$ and ${\mathcal H}_2$ are in
${\mathrm B}_n(\Omega)$ and $X:\cl H_1 \to \cl H_2$ is a bounded operator
intertwining these two operator tuples. Then $X$ must map the joint kernel of one
tuple in to the other, that is, $X K_1(\cdot, w) \x = K_2(\cdot,w) \Phi(w) \x$,
$\x\in \C^n,$ for some function $\Phi: \Omega \to \C^{n\times n}$. Assuming that the
kernel functions $K_1$ and $K_2$ are holomorphic in the first and anti-holomorphic
in the second variable, it follows, again as in \cite[pp. 472]{C-S}, that $\Phi$ is
anti-holomorphic. An easy calculation then shows that $X^*$ is the multiplication
operator $M_{\bar{\Phi}^{\rm tr}}$. If the two operator tuples are unitarily
equivalent then there exists an unitary operator $U$ intertwining them. Hence $U^*$
must be of the form $M_\Psi$ for some holomorphic function $\Psi$ such that
$\overline{\Psi(w)}^{\rm tr}$ maps the joint kernel of ${(\bm{M}-w)}^*$
isometrically onto the joint kernel of ${(\bm{M}-w)}^*$ for all $w\in \Omega$.  The
unitarity of $U$ is equivalent to the relation  $ K_1(\cdot,w)\x = U^* K_2(\cdot,w)
\overline{\Psi(w)}^{\rm tr}\x$ for all $w\in \Omega$ and $\x\in \C^n$.  It then
follows that
\begin{equation} \label{unitary rel kernel}
K_1(z,w) = \Psi(z) K_2(z,w)\overline{\Psi(w)}^{\rm tr},
\end{equation}
where $\Psi:\Omega_0 \subseteq \Omega \to {\mathcal G}{\mathcal L}(\F)$
is some holomorphic
function. Here, ${\mathcal G}{\mathcal L}(\F)$ denotes the group of all
invertible linear transformations on $\F$.

Conversely, if two kernels are related as above then the
corresponding tuples of multiplication \ops are unitarily
equivalent since
$$ M_i^*K(\cdot,w)\zeta=\bar{w_i}K(\cdot,w)\zeta,~~w\in
\Omega,~\zeta\in \F,
$$
where $(M_i f)(z)=z_if(z)$, $f\in {\mathcal H}$ for $1\leq i \leq m$.
\begin{rem} \label{projker}
We observe that if there is a self adjoint operator $X$ commuting
with the \dt $\bm{M}$ on the Hilbert space $\cl{H}$ then we must
have the relation $\overline{\Phi(z)}^{\rm tr} K(z,w) =
K(z,w) \Phi(w)$ for some anti-holomorphic function
$\Phi:\Omega \to \C^{n\times n}$. Hence if the kernel $K$ is
normalized then any projection $P$ commuting with the \dt $\bm{M}$
is induced by a constant function $\Phi$ such that $\Phi(0)$ is an
ordinary projection on $\C^n$.
\end{rem}

In conclusion, what is said above shows that a \dt of
operators in ${\mathrm B}_n(\Omega^*)$ admits a representation as the adjoint
of a \dt of multiplication
\ops on a \rk \hs of $\F$-valued holomorphic functions on $\Omega_0$,
where the \rk $K$ may be assumed to be normalized.  Conversely,
the adjoint of the \dt of multiplication \ops on the \rk \hs associated with a
normalized kernel $K$ on $\Omega$ belongs to
${\mathrm B}_n(\Omega^*)$ if certain additional conditions are imposed on
$K$ (cf. \cite{C-S}).

Our interest in the class ${\mathrm B}_n(\Omega)$ lies in the fact
that the Cowen-Douglas theorem \cite{C-D} provides a complete set of
unitary invariants for operators which belong to this class.  However,
these invariants are somewhat intractable.  Besides, often it is not
easy to verify that a given operator is in the class ${\mathrm
  B}_n(\Omega)$.  Although, we don't use the complete set of invariants
that \cite{C-D} provides, it is useful to ensure that the homogeneous
operators that arise from the jet construction are in this class.

\subsection{\sf The jet construction} \label{jet}
Let $\cl M$ be a Hilbert module over the algebra $\cl A(\Omega)$ for
$\Omega$ a bounded domain in $\C^m$.  Let $\cl M_k$ be the
submodule of functions in $\cl M$ vanishing to order $k>0$ on some analytic
hyper-surface $\cl Z$ in $\Omega$ -- the zero set of a holomorphic
function $\varphi$ in $\cl A(\Omega)$. A function $f$ on $\Omega$ is said to
vanish to order $k$ on $\cl Z$ if it can be written $f = \varphi^k g$ for some
holomorphic function $g$. The quotient module $\cl Q
= \cl M \ominus \cl M_k$ has been characterized in \cite{DMV}.
This was done by a generalization of the
approach in \cite{Aron} to allow vector-valued kernel Hilbert modules. The basic result in
\cite{DMV} is that $\cl Q$ can be characterized as such a vector-valued
kernel Hilbert space over the algebra $\cl A(\Omega)|_{\cl Z}$ of the restriction
of functions in $\cl A(\Omega)$ to $\cl Z$ and multiplication by $\varphi$ acts as
a nilpotent operator of order $k$.

For a fixed integer $k>0$, in this realization, $\cl M$ consists of $\C^k$-valued holomorphic functions,
and there is an $\C^{k\times k}$-valued
function $K(z,w)$ on $\Omega\times\Omega$ which is holomorphic
in $z$ and anti-holomorphic in $w$ such that
\begin{itemize}
\item[(1)] $K(\cdot,w)v$ is in $\cl M$ for $w$ in $\Omega$ and $v$ in $\bb C^k$;
\item[(2)] $\langle f,K(\cdot,w)v\rangle_{\cl M} = \langle f(w), v\rangle_{\mathbb
C^k}$ for $f$ in $\cl M$, $w$ in $\Omega$ and $v$ in $\bb C^k$; and \item[(3)] $\cl
A(\Omega) \cl M\subset \cl M$.
\end{itemize}

If we assume that $\cl{M}$ is in the class $\mathrm
B_1(\Omega)$, then it is possible to describe the quotient module via a
jet construction along the normal direction to the hypersurface
$\cl{Z}$.  The details are in \cite{DMV}.  In this approach, to
every positive definite kernel $K:\Omega\times \Omega \to \C$, we
associate a kernel $JK= \big (\!\!\big
(\partial_1^i\bar{\partial_1}^jK \big )\!\! \big )_{i,j=0}^{k-1}$,
where $\partial_1$ denotes differentiation along the normal direction
to $\cl{Z}$. Then we may equip
$$
J\cl{M} = \Big \{\mathbf f:= \sum_{i=0}^{k-1} \partial_1^i f \otimes
\varepsilon_i \in \mathcal M \otimes \C^k: f \in \cl{M} \Big \},
$$
where $\varepsilon_0, \ldots , \varepsilon_{k-1}$ are standard unit vectors in $\C^k$,
with a Hilbert space structure via the kernel $JK$. The module action is defined by
$\mathbf f \mapsto \mathbb J \mathbf f$ for $\mathbf f \in J\cl{M}$, where $ \mathbb J$ is the
array --
$$\mathbb J = \begin{pmatrix}
1 & \hdots & \hdots&\hdots & \hdots& 0\\
\partial_1 & 1 & & & & \vdots \\
\vdots & & \ddots & & & \vdots \\
\vdots &\binom{l}{j}\partial_1^{\ell - j} & &1& &\vdots\\
\vdots & & & &\ddots & 0 \\
\partial_1^{k-1}&\hdots &\hdots & \hdots & \hdots& 1 \cr
\end{pmatrix}
$$
with $0 \leq \ell,j \leq k-1$.
The module $J\cl{M}_{|{\rm res} ~\cl{Z}}$ which is the restriction of 
$J\cl{M}$ to $\cl{Z}$ is then shown to
be isomorphic to the quotient module $\cl{M} \ominus \mathcal M_k$.

We illustrate these results by means of an example.
Let  ${\cl M}^{(\alpha, \beta)}$ be the Hilbert module which corresponds to the reproducing kernel
$$B^{(\alpha, \beta)}(z, w)=\frac{1}{(1-z_1\bar{w}_1)^\alpha}\frac{1}{(1-z_2\bar{w}_2)^\beta},$$
$(z_1,z_2) \in \D^2$ and  $(w_1,w_2)\in \D^2$. Let ${\cl M}^{(\alpha, \beta)}_2$ be
the subspace of all functions in $\cl M^{(\alpha,\beta)}$ which vanish to order $2$
on the diagonal $\{(z,z)\,:\, z \in {\mathbb D}\} \subseteq {\mathbb D} \times
{\mathbb D}$. The quotient module $\cl Q :=\cl M^{(\alpha, \beta)} \ominus \cl
M^{(\alpha, \beta)}_2$ was described in \cite{DMTrans} using an orthonormal basis
for the quotient module $\cl Q$. This includes the calculation of the compression of
the two operators, $M_1:f\mapsto z_1f$ and $M_2: f \mapsto z_2f$ for $f\in \cl
M^{(\alpha, \beta)}$, on the quotient module $\cl Q$ (block weighted shift
operators) with respect to this orthonormal basis. These are homogeneous operators
in the class $\mathrm B_2(\D)$ which were first discovered by Wilkins \cite{wil}.

In  \cite{DMTrans}, an orthonormal basis $\left \{e_p^{(1)}, e_p^{(2)}\right
\}_{p=0}^\infty$ was constructed  in the quotient module $\cl M \ominus \cl
M^{(\alpha,\beta)}_2$. It was shown that the matrix
$$ M_p^{(1)} =
\begin{pmatrix}
\frac{\binom{-(\alpha+\beta)}{p}^{1/2}}{\binom{-(\alpha+\beta)}{p+1}^{1/2}} & 0\\
\big ( \frac{\beta}{\alpha}  \big )^{1/2} \frac{(\alpha + \beta
+1)^{1/2}}{\big ((\alpha + \beta +p) (\alpha+\beta + p + 1) \big
)^{1/2}} &
\frac{\binom{-(\alpha+\beta+2)}{p-1}^{1/2}}{\binom{-(\alpha+\beta+2)}{p}^{1/2}}
\end{pmatrix}$$
represents the operator $M_1$ which is multiplication by $z_1$
with respect to the orthonormal basis
$\{e_p^{(1)},e_p^{(2)}\}_{p=0}^\infty$. Similarly,
$$ M_p^{(2)} =
\begin{pmatrix}
\frac{\binom{-(\alpha+\beta)}{p}^{1/2}}{\binom{-(\alpha+\beta)}{p+1}^{1/2}} & 0\\
- \big ( \frac{\alpha}{\beta}  \big )^{1/2} \frac{(\alpha + \beta
+1)^{1/2}}{\big ((\alpha + \beta +p)(\alpha+\beta + p + 1) \big
)^{1/2}} &
\frac{\binom{-(\alpha+\beta+2)}{p-1}^{1/2}}{\binom{-(\alpha+\beta+2)}{p}^{1/2}}
\end{pmatrix}$$
represents the operator $M_2$ which is multiplication by $z_2$
with respect to the orthonormal basis
$\{e_p^{(1)},e_p^{(2)}\}_{p=0}^\infty$. Therefore, we  see that
$Q_1^{(p)} = \tfrac{1}{2}(M_1^{(p)} - M_2^{(p)})$ is a nilpotent
matrix of index $2$ while $Q_2^{(p)} =\tfrac{1}{2}(
M_1^{(p)}+M_2^{(p)})$ is a diagonal matrix in case $\beta=\alpha$.
These definitions naturally give a pair of operators $Q_1$ and
$Q_2$ on the quotient module $\cl Q^{(\alpha,\beta)}$. Let $f$ be a
function in the bi-disc algebra $\cl A(\D^2)$ and
$$f(u_1, u_2)= f_{0}(u_1) + f_{1}(u_1) u_2 + f_{2}(u_1) u_2^2 + \cdots
$$
be the Taylor expansion of the function $f$ with respect to the
coordinates $u_1 = \tfrac{z_1 + z_2}{2}$ and
$u_2 = \tfrac{z_1 - z_2}{2}$. Now the module
action for $f\in \cl A(\D^2)$ in the quotient module
$\cl Q^{(\alpha,\beta)}$ is then given by
\begin{eqnarray*}
f \cdot h &=& f(Q_1,Q_2) \cdot h \\
& = &  f_0(Q_1) \cdot h + f_1(Q_1) Q_2 \cdot h \\
&\stackrel{\rm def}{=}& \begin{pmatrix} f_0 & 0\\ f_1 & f_0
\end{pmatrix} \cdot
\begin{pmatrix} h_1 \\ h_2 \end{pmatrix},
\end{eqnarray*}
where $h= \tbinom{h_1}{h_2} \in \cl Q^{(\alpha,\beta)}$ is the
unique decomposition obtained from realizing the quotient module
as the direct sum $\cl Q^{(\alpha,\beta)} = \big (\cl M^{(\alpha,
\beta)} \ominus \cl M^{(\alpha,\beta)}_1 \big ) \oplus \big
(\cl M^{(\alpha,\beta)}_1\ominus \cl M^{(\alpha,\beta)}_2 \big )$,
where $\cl M^{(\alpha,\beta)}_i$, $i=1,2$, are the submodules in
$\cl M^{(\alpha,\beta)}$ consisting of all functions vanishing on
$\cl Z$ to order $1$ and $2$ respectively.

We now calculate the curvature $\cl K^{(\alpha, \beta)}$ for the
bundle $E^{(\alpha, \beta)}$ corresponding to the metric
$B^{(\alpha, \beta)}(\mathbf{u},\mathbf{u})$, where
$\mathbf{u}=(u_1,u_2) \in \D^2$.  The curvature $\cl K^{(\alpha,
\beta)}$ is easy to compute:
$$
\cl K^{(\alpha, \beta)} (u_1,u_2) = (1- |u_1 + u_2|^2)^{-2}
\begin{pmatrix} \alpha & \alpha \\ \alpha & \alpha
\end{pmatrix} +(1- |u_1 - u_2|^2)^{-2}
\begin{pmatrix} \beta & - \beta \\ -\beta & \beta
\end{pmatrix}.
$$
The restriction of the curvature to the hyper-surface $\{u_2=0\}$
is
$$
\cl K^{(\alpha, \beta)} (u_1,u_2)_{|u_2 = 0} = (1-|u_1|^2)^{-2}
\begin{pmatrix} \alpha + \beta & \alpha - \beta \\ \alpha - \beta &
\alpha + \beta
\end{pmatrix},
$$
where $u_1\in \D$. Thus we find that if $\alpha = \beta$, then the
curvature is of the form $2 \alpha (1-|u_1|^2)^{-2} I_2$.

We now describe the unitary map which is basic to the construction
of the quotient module, namely,
$$
h \mapsto  \sum_{\ell=0}^{k-1} \partial_1^{\ell} h\otimes
\varepsilon_{\ell} \Big |_{{\rm res ~}\triangle}
$$
for $h \in \cl M^{(\alpha, \beta)}$.  For $k=2$, it is enough to
describe this map just for the orthonormal basis
$\{e_p^{(1)}, e_p^{(2)}: p \geq 0\}$.  A simple calculation shows that
\begin{eqnarray}
e_p^{(1)}(z_1, z_2) &\mapsto& \begin{pmatrix}
\binom{-(\alpha+\beta)}{p}^{1/2} z_1^p \\
\beta \sqrt{\frac{p}{\alpha+\beta}} \binom{-(\alpha+\beta+1)}{p-1}^{1/2} z_1^{p-1}
\end{pmatrix} \nonumber \\
e_p^{(2)}(z_1, z_2) &\mapsto& \begin{pmatrix}
0\\
\sqrt{\frac{\alpha \beta}{\alpha+\beta}}\binom{-(\alpha+\beta+2)}{p-1}^{1/2} z_1^{p-1}
\end{pmatrix}.
\end{eqnarray}
This allows us to compute the $2\times 2$ matrix-valued kernel
function
$$
K_{\mathcal Q}({\mathbf z}, {\mathbf w}) =
\sum_{p=0}^\infty e_p^{(1)}({\mathbf z}) e_p^{(1)}({\mathbf w})^*
+ \sum_{p=0}^\infty e_p^{(2)}(\mathbf z) e_p^{(2)}(\mathbf w)^*,~{\mathbf z},
{\mathbf w}\in \D^2
$$
corresponding to the quotient module. Recall that $S(z,w):=(1-z\bar{w})^{-1}$ is the Szeg\"{o} kernel
for the unit disc $\mathbb D$.  We set $\bb S^r(z):= S(z,z)^r = (1-|z|^2)^{-r}$, $r>0$. 
A straight forward computation shows that
\begin{eqnarray*}
\lefteqn{K_{\mathcal Q}(\mathbf z, \mathbf z)_{|\mbox{res~ }\triangle} }\\
&=&\begin{pmatrix} \bb S(z)^{\alpha + \beta} &
\beta z \bb S(z)^{\alpha + \beta+1}\\
\beta \bar{z} \bb S(z)^{\alpha + \beta+1} &
\frac{\beta^2}{\alpha+\beta} \frac{d}{d |z|^2}
\big (|z|^2\bb S(z)^{\alpha + \beta+1}\big ) + \frac{\beta \alpha}{\alpha+\beta}
\bb S(z)^{\alpha + \beta+2}
\end{pmatrix} \\
&=& \big (\!\! \big ( \bb S(z_1)^\alpha {\partial^i}
\bar{\partial}^j
{\bb S(z_2)^\beta}_{|\mbox{res~ }\triangle} \big ) \!\! \big )_{i,j =0,1}\\
&=& (JK) (\mathbf z,\mathbf z)_{|\mbox{res~ }{\mathcal Z}},\:\: \mathbf z \in \D^2,
\end{eqnarray*}
where $\triangle=\{(z,z)\in \D^2 : z \in \D\}$. These calculations
give an explicit illustration of one of the main theorems on
quotient modules from  \cite[Theorem 3.4]{DMV}.

\section{Multiplier representations} \label{secmult}

Let $G$ be a locally compact second countable (lcsc) topological group
acting transitively on the domain $\Omega \subseteq
\bb{C}^m$. Let $\bb{C}^{n\times n}$ denote
the set of $n\times n$ matrices over the complex field $\bb{C}$.
We start with a cocycle $J$, that is, a holomorphic map $J_g:\Omega \to
\bb{C}^{n\times n}$ satisfying the cocycle relation
\begin{equation} \label{cocycle}
J_{gh}(z)= J_h(z)J_g(h\cdot z),~~\mbox{for all~} g,h\in G,~z\in\Omega,
\end{equation}
Let ${\rm Hol}(\Omega, \bb{C}^n)$ be the linear space consisting of all
holomorphic functions on $\Omega$ taking values in $\bb{C}^n$.
We then obtain a natural (left) action $U$ of the group $G$ on
${\rm Hol}(\Omega, \bb{C}^n)$:
\begin{equation}
  \label{eq:grp action}
  (U_{g^{-1}}f)(z) = J_g(z)f(g\cdot z),~f\in
  {\rm Hol}(\Omega,\bb{C}^n),~z\in \Omega.
\end{equation}
Let $e$ be the identity element of the group $G$.
Note that the cocycle condition \mbox{(\ref{cocycle})} implies, among other
things, $J_e(z)= J_e(z)^2$ for all $z\in \Omega$.

Let $\bb{K}\subseteq G$ be the compact subgroup which is the stabilizer of
$0$.  For $h,k$ in $\bb{K}$, we have $J_{kh}(0)=J_h (0)J_k (0)$ so that
$k\mapsto J_k (0)^{-1}$ is a \rep of $\bb{K}$ on $\bb{C}^n$.

A positive definite kernel $K$ on $\Omega$ defines an inner product
on some linear subspace of ${\rm Hol}(\Omega, \bb{C}^n)$.
The completion of this subspace is then a Hilbert space of holomorphic
functions on $\Omega$ (cf. \cite{Aron}). The natural action of
the group $G$ described above is seen to be unitary for an appropriate
choice of such a kernel.  Therefore, we first discuss these kernels in some
detail.

Let $\cl{H}$ be a functional Hilbert space consisting of holomorphic functions on
$\Omega$ possessing a \rk $K$.  We will always assume that the $m$ - tuple of
multiplication operators $\bm{M}=(M_1,\ldots, M_m)$ on the \hs $\cl{H}$ is bounded.
We also define the action of the group $G$ on the space of multiplication operators
-- $g\cdot M_f = M_{f\circ g}$ for $f \in \cl A(\Omega)$ and $g\in G$.  In
particular, we have $g\cdot \bm{M}=\bm{M}_g$. We will say that the \dt \bm{M} is
{\em G-homogeneous} if the operator $g\cdot \bm{M}$ is unitarily equivalent to
\bm{M} for all $g\in G$. $g\mapsto U_{g^{-1}}$ defined in \mbox{(\ref{eq:grp
action})} leaves $\cl{H}$ invariant. The following theorem says that the \rk of such
a \hs must be {\em quasi invariant} under the $G$ action.

A version of the following Theorem appears in \cite{KM} for the
unit disc. However, the proof here, which is taken from \cite{KM},
is for a more general domain $\Omega$ in $\C^m$.

\begin{thm} \label{main}
Suppose that $\mathcal H$ is a Hilbert  space which consists of
holomorphic functions on $\Omega$ and possesses a \rk $K$ on which
the $m$ - tuple $\bm{M}$ is irreducible and bounded.  Then the
following are equivalent.

\begin{enumerate}
\item The $m$ - tuple $\bm{M}$ is G-homogeneous.
\item The reproducing kernel $K$ of the \hs $\mathcal{H}$ transforms,
for some cocycle $J_g:\Omega \to \C^{n\times n}$, according to the rule
$$K(z,w) = J_g(z)K(g\cdot z, g\cdot w)J_g (w)^*,\:z,w\in\Omega.$$
\item The operator $U_{g^{-1}}: f \mapsto M_{J_g} f\circ g$ for $f\in \mathcal H$
is unitary.
\end{enumerate}

\end{thm}

\begin{proof}
Assuming that $K$  is quasi-invariant, that is, $K$ satisfies the transformation
rule, we see that the linear transformation $U$ defined in \mbox{(\ref{eq:grp action})}
is unitary. To prove this, note that
\begin{align*}
\inner{U_{g^{-1}}K(z, w)\x}{U_{g^{-1}}K(z, w^\prime)\y} &=
\inner{J_g(z) K(g\cdot z, w)\x}{J_g(z)K(g\cdot z, w^\prime)\y}\\
&= \inner{K(z,\tilde{w}){{J_g(\tilde{w})}^*}^{-1}\x}
{K( z,\tilde{w}^\prime) {J_g(\tilde{w}^\prime)^*}^{-1}\y}\\
&= \inner{K(\tilde{w}^\prime,
\tilde{w}){{J_g(\tilde{w})}^*}^{-1}\x}
{ {J_g(\tilde{w}^\prime)^*}^{-1}\y}\\
&=\inner{{J_g(\tilde{w}^\prime)}^{-1}K(\tilde{w}^\prime,
\tilde{w}){{J_g(\tilde{w})}^*}^{-1}\x}{\y}\\
&=\inner{K(g\cdot\tilde{w}^\prime ,
g\cdot\tilde{w})\x}{\y},
\end{align*}
where $\tilde{w} = g^{-1}\cdot w$ and $\tilde{w}^\prime = g^{-1}\cdot{w}^\prime$.
Hence
$$\inner{K(g\cdot\tilde{w}^\prime , g\cdot\tilde{w})\x}{\y}
= \inner{K(w^\prime, w)\x}{\y}.
$$
It follows that the map $U_{g^{-1}}$ is isometric. On the other
hand, if $U$ of \mbox{(\ref{eq:grp action})} is unitary then the
reproducing kernel $K$ of the  \hs $\mathcal{H}$ satisfies
\begin{equation}
  \label{eq:trans rule}
 K(z,w) = J_g(z)K(g\cdot z, g\cdot w)J_g(w)^*.
\end{equation}
This follows from the fact that the reproducing kernel has the
expansion \mbox{(\ref{sumonb})}
for some orthonormal basis  $\{e_\ell:\ell \geq 0\}$ in $\cl{H}$.
The uniqueness of the \rk implies that the expansion is
independent of the choice of the orthonormal basis.  Consequently,
we also have $K(z,w) = \sum_{\ell=0} (U_{g^{-1}}e_\ell)
(z)(U_{g^{-1}}e_\ell) (w)^*$ which verifies the equation
\mbox{(\ref{eq:trans rule})}. Thus we have shown that $U$ is
unitary if and only if the \rk $K$ transforms according to
\mbox{(\ref{eq:trans rule})}.

We now show that the \dt \bm{M} is homogeneous if and only if $f\mapsto M_{J_g}
f\circ g$ is unitary. The eigenvector at $w$ for $g.\bm{M}$ is clearly $K(\cdot,
g^{-1}\cdot w)$. It is not hard, using the unitary operator $U_\Gamma$ in
\mbox{(\ref{general construction})}, to see that that $g^{-1} \cdot \bm{M}$ is
unitarily equivalent to \bm{M} on a Hilbert space $\cl H_g$ whose \rk is $K_g(z,w) =
K(g\cdot z, g \cdot w)$ and the unitary $U_\Gamma$ is given by $f \mapsto f \circ g$
for $f \in \cl H$.  However, the homogeneity of the \dt \bm{M} is equivalent to the
existence of a unitary operator intertwining the \dt of multiplication on the two
Hilbert spaces $\cl H$ and $\cl H_g$. As we have pointed out in section \ref{CD},
this unitary operator is induced by a multiplication operator $M_{J_g}$, where $J_g$
is a holomorphic function (depends on $g$) such that
$K_g(z,w)=J_g(z)K(z,w)\overline{J_g(w)}^{\rm tr}$. The composition of these two
unitaries is  $f \mapsto M_{J_g} f \circ g$ and is therefore a unitary.
\end{proof}

The discussion below and the Corollary following it is implicit in \cite{KM}.
Let $g_z$ be an element of $G$ which maps $0$ to $z$, that is $g_z
\cdot 0 = z$. We could then try to define possible kernel functions
$K:\Omega\times \Omega\to \bb{C}^{n\times n}$ satisfying the
transformation rule \mbox{(\ref{eq:trans rule})} via the requirement
\begin{equation}
  \label{eq:homker def}
  K(g_z\cdot 0, g_z\cdot 0) = (J_{g_z}(0))^{-1}K(0,0)(J_{g_z}(0)^*)^{-1},
\end{equation}
choosing any positive \op $K(0,0)$ on $\bb{C}^n$ which commutes with $J_k(0)$
for all $k\in \bb{K}$. Then the equation \mbox{(\ref{eq:homker def})}
determines the function $K$ unambiguously as long as $J_k(0)$ is
unitary for $k\in \bb{K}$.  Pick $g\in G$ such that $g\cdot 0=z$.
Then $g=g_z k$ for some $k\in \bb{K}$. Hence
\begin{eqnarray*}
  K(g_z k\cdot 0, g_z k\cdot 0)&=& (J_{g_zk} (0))^{-1}K(0,0)(J_{g_zk}
  (0)^*)^{-1}\\
  &=& \big ( J_k(0)J_{g_z}(k\cdot 0) \big )^{-1}K(0,0)\big (
  J_{g_z}(k\cdot 0)^*  J_k(0)^*\big )^{-1} \\
  &=& (J_{g_z}(0))^{-1} (J_k(0))^{-1} K(0,0)
  (J_k(0)^*)^{-1}(J_{g_z}(0)^*)^{-1}\\
  &=& (J_{g_z}(0))^{-1} K(0,0)(J_{g_z}(0)^*)^{-1}\\
  &=& K(g_z\cdot 0, g_z\cdot 0)
\end{eqnarray*}

Given the definition \mbox{(\ref{eq:homker def})}, where the choice of
$K(0,0)=A$ involves as many parameters as the number of irreducible
representations of the form $k \mapsto J_k(0)^{-1}$ of the compact
group $\bb{K}$, one can polarize (\mbox{\ref{eq:homker def})} to get
$K(z,w)$.  In this approach, one has to find a way of
determining if $K$ is non-negative definite, or for that matter,
if $K(\cdot,w)$ is holomorphic on all of $\Omega$ for each fixed but
arbitrary $w\in \Omega$.   However, it is evident from the definition
\mbox{(\ref{eq:homker def})} that
\begin{eqnarray*}
  K(h\cdot z,h\cdot z)&=&J_h(g_z\cdot 0)^{-1}J_{g_z}(0)^{-1}A
  {J_{g_z}(0)^*}^{-1}(J_h(g_z\cdot 0)^*)^{-1}\\
  &=&J_h(z)^{-1}K(z,z){J_h(z)^*}^{-1}
\end{eqnarray*}
for all $h\in G$.  Polarizing this equality, we obtain
$$K(h\cdot z,h\cdot w)=J_h(z)^{-1} K(z,w){J_h(w)^*}^{-1}$$ which is the identity
\mbox{(\ref{eq:trans rule})}.  It is also clear that the linear span of the
set $\{K(\cdot,w)\zeta:~ w\in \Omega,~\zeta\in \bb{C}^n\}$ is stable
under the action \mbox{(\ref{eq:grp action})} of $G$:
$$
g\mapsto J_g(z)K(g\cdot z,w)\zeta = K(z, g^{-1}\cdot w){J_g(g^{-1}w)^*}^{-1}\zeta,
$$
where $J_g(g^{-1}w)^{*-1}\zeta$ is a fixed element of $\bb{C}^n$.
\begin{cor} \label{}
If $J:G\times \Omega\to\C^{n\times n}$ is a cocycle and $g_z$ is an element of
$G$ which maps $0$ to $z$ then the kernel $K : \Omega \times \Omega
\to \C^{n\times n}$ defined by the requirement
$$
K(g_z\cdot 0, g_z\cdot 0) = (J_{g_z}(0))^{-1}K(0,0)(J_{g_z}(0)^*)^{-1}
$$
is quasi-invariant, that is, it transforms according to \mbox{(\ref{eq:trans rule})}.
\end{cor}

\section{Irreducibility}
In the section \ref{CD}, we have already pointed out that
any Hilbert space $\mathcal H$ of scalar valued holomorphic
functions on $\Omega\subset{\mathbb C}^m$ with a reproducing
kernel $K$ determines a line bundle $\mathcal E$ on
$\Omega^*:=\{\bar w:w\in \Omega\}$. The fibre of $E$ at $\bar w\in
\Omega^*$ is spanned by $K(.,w).$ We can now construct a rank
$(n+1)$ vector bundle $J^{(n+1)}\mathcal E$ over $\Omega^*$.  A
holomorphic frame for this bundle is $\{{\bar\partial}^l_2 K(.,w):0\leq
l\leq k,w\in\Omega\},$ and as usual, this frame determines a
metric for the bundle which we denote by $J^{(n+1)}K,$ where
$$J^{(n+1)}K(w,w)=\big
( \!\! \big (\langle \bar\partial_2^j K(., w), \bar\partial_2^i
K(., w) \rangle \big )\!\! \big )_{i,j=0}^n=\big ( \!\! \big
(\bar\partial_2^j\partial_2^iK(w , w)  \big )\!\! \big
)_{i,j=0}^n,w\in \Omega.$$

Recall that the kernel function on
$\D^2$, $B^{(\alpha,\beta)}:\D^2\times\D^2\longrightarrow\mathbb C$ is defined by
$$B^{(\alpha,\beta)}(z,w)=(1-z_1\bar{w_1})^{-\alpha}(1-z_2\bar{w_2})^{-\beta},$$
for $z=(z_1, z_2)\in \D^2$ and $w=(w_1,w_2)\in \D^2$, $\alpha,
\beta > 0.$ Take $\Omega=\D^2, K= B^{(\alpha,\beta)}$. Notice that
the Hilbert space $\cl M^{(\alpha,\beta)}$ corresponding to the
kernel function $B^{(\alpha,\beta)}$ is the tensor product of the
two Hilbert spaces $\cl M^{(\alpha)}$ and $\cl M^{(\beta)}$. These
are determined by the two kernel functions
$B^{(\alpha)}(z,w)=(1-z\bar{w})^{-\alpha}$ and
$B^{(\beta)}(z,w)=(1-z\bar{w})^{-\beta}$, $z,w\in \D$,
respectively.

It follows from \cite{DMV} that
$h_{n+1}(z)=J^{(n+1)}B^{(\alpha,\beta)}(z,z)_{|{\rm res}~\triangle}$ is a metric for the 
Hermitian anti-holomorphic vector bundle $J^{(n+1)}\mathcal E_{|{\rm res}~\triangle}$ over 
$\triangle=\{(z,z): z\in\D\}\subseteq\D^2.$ However, $J^{(n+1)}\mathcal E_{|{\rm res}~\triangle}$ is a
Hermitian holomorphic vector bundle over $\triangle^*=\{(\bar z,\bar z):z\in\D\}$, that is, 
$\bar z$ is the holomorphic variable in this description. Thus $\partial f=0$ if and only if
$f$ is holomorphic on $\triangle^*$. To restore the usual meaning of $\partial$ and
$\bar\partial,$ we interchange the roles of $z$ and $\bar z$ in the metric  which
amounts to replacing $h_{n+1}$ by its transpose.

As shown in \cite{DMV}, this Hermitian anti-holomorphic vector bundle
$J^{(n+1)}\mathcal E_{|{\rm res}~\triangle}$ defined over the diagonal subset $\triangle$  
of  the bidisc $\bb D^2$ gives rise to a reproducing kernel Hilbert space $J^{(n+1)}\mathcal H$.   
The  reproducing kernel for this Hilbert space is $J^{(n+1)}B^{(\alpha,\beta)}(z,w)$ which is 
obtained by polarizing $J^{(n+1)}B^{(\alpha,\beta)}(z,z)=h_{n+1}(z)^t$.

\begin{lem} \label{newpair} Let $\alpha, \beta$ be two positive
 real numbers and $n\geq 1$ be an integer. Let $\cl M_n$ be the
ortho-complement of the subspace of
 $\cl M^{(\alpha)} \otimes \cl M^{(\beta)}$ (viewed as a Hilbert space of
 analytic functions on the bi-disc $\D \times \D$) consisting of all
 the functions vanishing to order $k$ on the diagonally embedded unit
 disc $\triangle :=\{(z,z):z\in \D \}$. The
compressions to $\cl M_n$ of $M^{(\alpha)} \otimes I$ and
$I \otimes M^{(\beta)}$ are
 homogeneous operators with a common associated representation.
\end{lem}

\begin{proof} For each real number $\alpha > 0$, let $\cl M^{(\a)}$
be the Hilbert space completion of the inner product space spanned by
$\{ f_k : k \in \mathbb Z^+\}$ where  the $f_k$'s are mutually orthogonal
vectors with norms given by
$$\|f_k\|^2=\frac{\Gamma(1+k)}{\Gamma(\a + k)},\: k \in \mathbb Z^+.$$
(Upto scaling of the norm, this Hilbert space may be identified,
via non-tangential boundary values, with the Hilbert space of
analytic functions on $\D$ with reproducing kernel $(z,w)\mapsto
(1-z\bar{w})^{-\alpha}$.) The representation $D_\alpha^+$ lives on
$\cl M^{(\alpha)}$, and is given (at least on the linear span of the
$f_k$'s) by the formula
$$D_{\alpha}^+(\varphi^{-1})f = (\varphi^\prime )^{\alpha/2}f\circ
\varphi,\: \varphi \in \mbox{M\"{o}b}.$$

Clearly the subspace $\cl M_n$ is invariant under the Discrete series representation
$\pi := D_\alpha^+ \otimes D_\beta^+$ associated with both the operators
$M^{(\alpha)} \otimes I$ and $I\otimes M^{(\beta)}$.  It is also co-invariant under
these two operators. An application of Proposition 2.4 in \cite{const} completes 
the proof of the lemma.
\end{proof}

The subspace $\cl M_n$ consists of those functions $f \in \cl M$ which vanish 
on $\triangle$ along with their first $n$ derivatives with respect to $z_2$. 
As it turns out, the compressions to $\cl M\ominus \cl M_n$ of $M^{(\alpha)} \otimes I$
is the multiplication operator on the Hilbert space ${J^{(n+1)}\mathcal H}_{|{\rm res}~\triangle}$
which we denote $M^{(\a,\b)}$.  An application of \cite[Proposition 3.6]{DMV} shows that the 
adjoint $M^*$ of the multiplication operator $M$ is in $\mathrm B_{n+1}(\D)$. 

\begin{thm}\label{tirr}
The  multiplication operator $M:=M^{(\a,\b)}$ is irreducible. 
\end{thm} 

Th proof of this theorem will be facilitated by a series of lemmas which are proved in the sequel.  
We set, for now,  $K(z,w)=J^{(n+1)}B^{(\alpha,\beta)}(z,w)$.  Let $\widetilde
K(z,w)=K(0,0)^{-1/2}K(z,w)K(0,0)^{-1/2},$ so that $\widetilde K(0,0)=I.$  
Also, let  $\widetilde{\widetilde K}(z,w)=\widetilde K(z,0)^{-1}\widetilde
K(z,w)\widetilde K(0,w)^{-1}$.  This ensures that $\widetilde{\widetilde K}(z,0)=I$ for
$z\in\D,$ that is, $\widetilde{\widetilde K}$ is a normalized kernel. 
Each of the kernels $K$, $\widetilde{K} \mbox{~and~} \widetilde{\widetilde K}$ admit 
a power series expansion, say,  $K(z,w)=\sum_{m,\,p\,\geq \,0}a_{mp}\,z^m{\bar w}^p$, 
$\widetilde K(z,w)=\sum_{m,\,p\,\geq \,0}{\widetilde a}_{mp}\,z^m{\bar w}^p,$  and 
$\widetilde{\widetilde K}(z,w)=\sum_{m,\,p\,\geq \,0}\widetilde{\widetilde a}_{mp}\,z^m{\bar
w}^p $ for $z,\,w\in\D$, respectively.  Here the coefficients $a_{mp}$ and $\widetilde
a_{mp}$ and $\widetilde{\widetilde a}_{mp}$ are $(n+1)\times(n+1)$ matrices 
for $m,\,p \geq 0$.  In particular,  $\widetilde a_{mp} = K(0,0)^{-1/2}a_{mp}\,K(0,0)^{-1/2} =  a_{00}^{-1/2}a_{mp}\,a_{00}^{-1/2}$ for $m,p\geq 0$.
Also, let us write $K(z,w)^{-1}=\sum_{m,\,p\,\geq \,0}{ b}_{mp}\,z^m{\bar w}^p$ and $\widetilde
K(z,w)^{-1}=\sum_{m,\,p\,\geq\, 0}{ \widetilde b}_{mp}\,z^m{\bar w}^p$, $z, w\in \D.$  Again, 
the coefficients $b_{mp}$ and $\widetilde b_{mp}$ are $(n+1)\times(n+1)$
matrices for $m,\,p \geq 0$.  However, $\wi{\wi a}_{00}=I$ and $\wi{\wi a}_{m0}=\wi{\wi a}_{0p}=0$
for $m,\,p\geq 1$.

The following lemma is from \cite[Theorem 3.7, Remark 3.8 and Lemma
3.9]{C-S}. The proof was discussed in section \ref {CD}.
\begin{lem} \label{kur}
The multiplication operators on Hilbert spaces $\mathcal H_1$ and
$\mathcal H_2$ with reproducing kernels $K_1(z,w)$ and $K_2(z,w)$
respectively, are \u\ if and only if
$K_2(z,w)=\Psi(z)K_1(z,w)\overline{\Psi(w)}^t$, where $\Psi$ is an
invertible matrix-valued holomorphic function. 
\end{lem}
The proof of the  lemma below appears in \cite[Lemma 5.2]{KM}
and is discussed in section \ref{CD}, see Remark \ref{projker}.
\begin{lem}
The multiplication operator $ M$ on the Hilbert space $ {\mathcal H}$ with
reproducing kernel $K$ is irreducible if and only if there is no non-trivial
projection $P$ on $\C^{n+1}$ commuting with all the coefficients in the power series
expansion of the  normalized kernel $\widetilde{\widetilde K}(z,w)$. 
\end{lem}

We will prove irreducibility of $ M$ by showing that only operators on $\C^{n+1}$
which commutes with all the coefficients of $\widetilde {\wi K}(z,w)$ are scalars.
It turns out that the coefficients of $z^k\bar w$ for $2\leq k\leq n+1,$ that is,
the coefficients $\widetilde{\widetilde a}_{k1}$ for $2\leq k\leq n+1$ are
sufficient to reach the desired conclusion.
\begin{lem}\label{coeff}
The coefficient of $z^k\bar w$ is  $\widetilde{\widetilde
a}_{k1}=\d_{s=1}^k\widetilde b_{s0}\widetilde a_{k-s,1}+\widetilde a_{k1}$ for
$1\leq k\leq n+1$.
\end{lem}
\begin{proof}
Let us denote the  coefficient of $z^k{\bar w}^l$ in the power series expansion of
$\wi{\wi K}(z,w)$ is $\wi{\wi a}_{kl} $ for $k,l\geq 0.$  We see that 
\Bea \lefteqn{
{\wi{\wi a}_{kl} = \d_{s=0}^k}\d_{t=0}^l \wi b_{s0}\wi
a_{k-s,l-t}\wi b_{0t}}\\
&&=\d_{s=1}^k\d_{t=1}^l\wi a_{s0}\wi a_{k-s,l-t}\wi b_{0t}
+\d_{s=1}^k\wi b_{s0}\wi a_{k-s,l}+\d_{t=1}^l\wi a_{k,l-t}\wi
b_{0t}+\wi a_{kl}
\Eea 
as $\widetilde a_{00}=\widetilde b_{00}=I.$  Also, 
\Bea \lefteqn{\widetilde{\widetilde a}_{k1}
=\sum_{s=1}^k\widetilde b_{s0}\widetilde a_{k-s,0}\widetilde
b_{01}+\d_{s=1}^k\widetilde b_{s0}\widetilde a_{k-s,1}+\widetilde
a_{k0}\widetilde b_{01}+\widetilde a_{k1}}\\
&&=\big(\sum_{s=0}^k\widetilde b_{s0}\widetilde
a_{k-s,0}\big)\widetilde b_{01}+\sum_{s=1}^k\widetilde
b_{s0}\widetilde
a_{k-s,1}+\widetilde a_{k1}\\
&&=\sum_{s=1}^k\widetilde b_{s0}\widetilde a_{k-s,1}+\widetilde
a_{k1} 
\Eea 
as $\widetilde b_{00}=I$ and coefficient of $z^k $ in
$\widetilde K(z,w)^{-1}\widetilde K(z,w)=\sum_{s=0}^k\widetilde
b_{s0}\widetilde a_{k-s,0}=0$ for $k\geq 1.$
\end{proof}

Now we compute some of the coefficients of $K(z,w)$ which are useful in computing
$\widetilde{\widetilde a}_{k1}$.   In what follows, we will compute only the
non-zero entries of the matrices involved, that is, {\em all those entries which are not specified
are assumed to be zero.}  

\begin{lem}\label{relcoeff}
$(a_{00})_{kk}=k!(\b)_{k}$ for $0\leq k\leq n,$
$(a_{m0})_{r,r+m}=\frac{(m+r)!}{m!}(\b)_{m+r}$ and
$(a_{m+1,1})_{r,r+m}=\frac{(m+r)!}{m!}(\b)_{m+r}\big(\a+(1+\frac{r}{m+1})(\b+m+r)\big)$
for $0\leq r\leq n-m,0\leq m\leq n,$  where
$(x)_0=1,(x)_d=x(x+1)\ldots (x+d-1)$, for any positive integer $d$,
is the Pochhammer symbol.
\end{lem}
\begin{proof} The coefficient of $z^p\bar w^q$ in $J^{(n+1)}B^{(\alpha,\beta)}(z,w)$ 
is the same as the coefficient of $z^p\bar z^q$ in $J^{(n+1)}B^{(\alpha,\beta)}(z,z)$. So,
$(a_{00})_{kk} =$ constant term in $\bar\partial_2^{k}\partial_2^k  
\big(\bb S(z_1)^\a\bb S(z_2)^\b\big)|_{\triangle}.$  Now,  
\Bea
\lefteqn{\bar\partial_2^{k}\partial_2^k\big(\bb S(z_1)^\a\bb S(z_2)^\b\big)|_{\triangle} = \bar\partial_2^{k}\partial_2^k\big(\bb S(z_1)^\a\bb S(z_2)^\b\big)|_{\triangle}}\\
&&=\bb S(z_1)^\a(\b)_k\bar\partial_2^{k}\big(\bb S(z_2)^{\b+k}\bar
z_2^k\big)|_{\triangle}\\
&&=\big(\bb S(z_1)^\a(\b)_{k}\sum_{l=0}^k\binom{k}{l}
\bar\partial_2^{k-l}(\bb S(z_2)^{\b+k})\bar\partial_2^{l}(\bar
{z_2}^k)\big)|_{\triangle}\\ &&= \big(\bb S(z_1)^\a(\b)_k\d_{l=0}^k
\binom{k}{l}(\b+k)_{k-l}\bb S(z_2)^{\b+k+(k-l)}{z_2}^{k-l}
l!\binom{k}{l}\bar{z_2}^{k-l}\big)|_{\triangle},
\Eea 
that is, $\big(a_{00}\big)_{kk}= k!(\b)_{k}$ for $0\leq k\leq n.$

We see that the coefficient of $z^m$ \mbox{~in~}
$\bar\partial_2^{m+r}\partial_2^r\big(\bb S(z_1)^\a\bb S(z_2)^\b\big)|_{\triangle}$
is $(a_{m0})_{r,r+m}$. Thus
\Bea 
\lefteqn{\bar\partial_2^{m+r}\partial_2^r\big(\bb S(z_1)^\a\bb S(z_2)^\b\big)|_{\triangle}
=\bb S(z_1)^\a(\b)_r\bar\partial_2^{m+r}\big(\bb S(z_2)^{\b+r}\bar
z_2^r\big)|_{\triangle}}\\
&&= \big(\bb S(z_1)^\a(\b)_r\d_{l=0}^{m+r}\tbinom{m+r}{l}
\bar\partial_2^{m+r-l}(\bb S(z_2)^{\b+r})\bar\partial_2^{l}(\bar
{z_2}^r)\big)|_{\triangle}\\
&&=
\big(\bb S(z_1)^\a(\b)_r\sum_{l=0}^{m+r}\tbinom{m+r}{l}(\b+r)_{m+r-l}
\bb S(z_2)^{\b+2r+m-l)}{z_2}^{m+r-l}
l!\binom{r}{l}\bar{z_2}^{r-l}\big)|_{\triangle}.\\
\Eea 
Therefore, the term containing $z^m$ occurs only when $l=r$ in the sum above, that is,
$(a_{m0})_{r,r+m}=(\b)_r\binom{m+r}{r}(\b+r)_{m}r!=\frac{(m+r)!}{m!}(\b)_{m+r},$ 
for $0\leq r\leq n-m, 0\leq m\leq n.$

Coefficient of $z^{m+1}\bar z$ in
$\bar\partial_2^{m+r}\partial_2^r\big(\bb S(z_1)^\a\bb S(z_2)^\b\big)|_{\triangle}
$ is  $(a_{m+1,1})_{r,r+m}.$ For any real analytic function
$f$ on $\D,$ for now, let $\big (f(z,\bar z)\big )_{(p,q)}$ denote the coefficient of $z^p{\bar z}^q$ in
$f(z,\bar z)$.  We have 

\Bea \lefteqn{\big(a_{m+1,1}\big)_{r,r+m}
=\big(\bar\partial_2^{m+r}\partial_2^r\big(\bb S(z_1)^\a\bb S(z_2)^\b\big)|_{\triangle}\big)_{(m+1,1)}}\\
&=&\Big((\b)_r\sum_{l=0}^{m+r}\tbinom{m+r}{l}(\b+r)_{m+r-l}\bb S(z)^{\a+\b+r+(m+r-l)}{z}^{m+r-l}
l!\binom{r}{l}\bar{z}^{r-l}\Big)_{(m+1,1)}
\Eea 
The terms containing $z^{m+1}\bar z$ occurs in the sum above, only when $l=r$ and $l=r-1$, that is,
\Bea
\lefteqn{\big(a_{m+1,1}\big)_{r,r+m}
=\big((\b)_r r!\bigg(\tbinom{m+r}{r}(\b+r)_m\bb S(z)^{\a+\b+m+r}z^m}\\
&& \hspace*{80pt}+\tbinom{m+r}{r-1}(\b+r)_{m+1}\bb S(z)^{\a+\b+m+r+1}z^{m+1}\bar
z\bigg)\big)_{(m+1,1)}\\
&=&\big((\b)_rr!\bigg(\frac{(m+r)!}{r!m!}(\b+r)_{m}(1+(\a+\b+m+r)|z|^2)z^m\\
&&\hspace*{90pt}+\frac{(m+r)!r}{r!(m+1)!}(\b+r)_{m+1}\bb S(z)^{\a+\b+m+r+1}z^{m+1}\bar
z\bigg)\big)_{(m+1,1)}\\
&=&\frac{(m+r)!}{m!}(\b)_{m+r}\big((\a+\b+m+r)+\frac{r}{m+1}(\b+m+r)\big)\\
&=& \frac{(m+r)!}{m!}(\b)_{m+r}\big(\a+(1+\frac{r}{m+1})(\b+m+r)\big),
\Eea 
for $0\leq r\leq n-m, 0\leq m\leq n,$ where we have followed the convention: $\binom{p}{q}=0$
for a negative integer $q.$
\end{proof}

\begin{lem}\label{impcoeff} Let $c_{k0}$ denote $a_{00}^{1/2}\widetilde b_{k0}a_{00}^{1/2}$.  
For $0\leq r\leq n-k,\,0\leq k\leq n$,  $(c_{k0})_{r,r+k}=\frac{(-1)^k(r+k)!}{k!}(\b)_{r+k}$.
\end{lem}
\begin{proof}
We have  
$\widetilde K(z,w)^{-1}=a_{00}^{1/2}K(z,w)^{-1}a_{00}^{1/2} = \sum_{mn\geq 0} 
\big ( a_{00}^{1/2} \widetilde b_{mn} a_{00}^{1/2} \big ) z^m{\bar w}^n.$  Hence 
$\widetilde b_{mn}= a_{00}^{1/2}b_{mn}a_{00}^{1/2}$ for $m,n\geq 0$. By
invertibility of $a_{00},$ we see that $\wi b_{k0}$ and $c_{k0}$ uniquely
determine each other for $k\geq 0$. Since $(\widetilde b_{k0})_{k\geq
0}$ are uniquely determined as the coefficients of power series
expansion of $\widetilde K(z,w)^{-1}$, it is enough to prove that
$\d_{l=0} ^m\widetilde a_{m-l,0}\widetilde b_{l0}=0$ for $1\leq
m\leq n.$ Equivalently, we must show that 
$\d_{l=0}^m(a_{00}^{-1/2}a_{m-l,0}a_{00}^{-1/2})(a_{00}^{-1/2}
c_{l0}a_{00}^{-1/2})=0$ which amounts to showing $a_{00}^{-1/2}\big(\d_{l=0}^m
a_{m-l,0}a_{00}^{-1} c_{l0}\big)a_{00}^{-1/2}=0$ for $1\leq m\leq
n.$ It follows from Lemma \ref{relcoeff} that
$(a_{m-l,0})_{r,r+(m-l)}=\frac{(m-l+r)!}{(m-l)!}(\b)_{m-l+r}$ and
$(a_{00})_{rr}=r!(\b)_{r}.$  Therefore
\Bea
\lefteqn{(a_{m-l,0}a_{00}^{-1})_{r,r+(m-l)}=(a_{m-l,0})_{r,r+(m-l)}
(a_{00}^{-1})_{r+(m-l),r+(m-l)}}\\
&&\hspace*{75pt}=\frac{(m-l+r)!}{(m-l)!}
(\b)_{m-l+r}\big((m-l+r)!(\b)_{m-l+r}\big)^{-1}\\
&&\hspace*{75pt}=\frac{1}{(m-l)!}.
\Eea
We also have
\Bea
\lefteqn{(a_{m-l,0}a_{00}^{-1}c_{l0})_{r,r+m}=(a_{m-l,0}a_{00}^{-1})_{r,r+(m-l)}
(c_{l0})_{r+(m-l),r+(m-l)+l}}\\
&&\hspace*{71pt}=\frac{1}{(m-l)!}\frac{(-1)^l(r+m)!}{l!}(\b)_{r+m}\\
&&\hspace*{71pt}=\frac{(-1)^l(r+m)!}{(m-l)!l!}(\b)_{r+m}
\Eea 
for $0\leq l\leq m,0\leq r\leq n-m,1\leq m\leq n$.  Now observe that
\Bea
(\d_{l=0}^m
a_{m-l,0}a_{00}^{-1}c_{l0})_{r,r+m}&=& (r+m)!(\b)_{m+r}\d_{l=0}^m
\frac{(-1)^l}{(m-l)!l!}\\
&=&\frac{(r+m)!}{m!}(\b)_{m+r}\d_{l=0}^m(-1)^l\binom{m}{l}\\
&=&0,
\Eea
which completes the proof of this lemma.
 \end{proof}
 \begin{lem}\label{oncoeff}
 $(\widetilde{\widetilde a}_{k1})_{n-k+1,n}$ is a non-zero
 real number, for $2\leq k\leq n+1,n\geq 1$. All other entries of $~~\widetilde{\widetilde
 a}_{k1}$ are zero.
 \end{lem}
 \begin{proof}
 From Lemma \ref{coeff} and Lemma \ref{impcoeff}, we know that 
 \Bea
 \widetilde{\widetilde
a}_{k1} &=& \d_{s=1}^k\widetilde b_{s0}\widetilde a_{k-s,1}+\widetilde a_{k1}\\
& = &\d_{s=1}^k(a_{00}^{-1/2}c_{s0}a_{00}^{-1/2})(a_{00}^{-1/2}a_{k-s,1}a_{00}^{-1/2})
+a_{00}^{-1/2}a_{k1}a_{00}^{-1/2}.
\Eea 
Consequently,  $a_{00}^{1/2}\widetilde{\widetilde a}_{k1}a_{00}^{1/2}=\d_{s=1}^k
c_{s0}a_{00}^{-1}a_{k-s,1} +a_{k1}$ for $1\leq k\leq n+1.$  By
Lemma \ref{relcoeff} and Lemma \ref{impcoeff}, we have
\Bea
(c_{s0}a_{00}^{-1})_{r,r+s} &=&(c_{s0})_{r,r+s}(a_{00}^{-1})_{r+s,r+s} \\
&=&\frac{(-1)^s(r+s)!}{s!}(\b)_{r+s}\big((r+s)!(\b)_{r+s}\big)^{-1}\\
&=&\frac{(-1)^s}{s!},
\Eea 
for $0\leq r\leq n-s, 0\leq s\leq k, 1\leq k\leq n+1.$
\Bea
\lefteqn{(a_{k-s,1})_{r,r+(k-s-1)}=}\\ 
&&\frac{(k+r-s-1)!}{(k-s-1)!}(\b)_{r+k-s-1} \big(\a+(1+\frac{r}{k-s})(\b+r+k-s-1)\big),
\Eea 
for $k-s-1\geq 0, 2\leq k\leq n+1.$ Now,
\Bea
\lefteqn{(c_{s0}a_{00}^{-1}a_{k-s,1})_{r+s,r+s+(k-s-1)}}\\
&=&(c_{s0}a_{00}^{-1})_{r,r+s}(a_{k-s,1})_{r+s,r+s+(k-s-1)}\\
&=& \frac{(-1)^s}{s!}\frac{(r+k-1)!}{(k-s-1)!}(\b)_{r+k-1}
\big(\a+(1+\frac{r+s}{k-s})(\b+r+k-1)\big),
\Eea 
for $1\leq s\leq k-1,0\leq r\leq n-k+1,1\leq k\leq n+1.$  Hence 
\Bea
\lefteqn{(c_{s0}a_{00}^{-1}a_{k-s,1})_{r+s,r+k-1}=}\\
&&\frac{(-1)^s}{s!}\frac{(r+k-1)!}{(k-s-1)!}(\b)_{r+k-1}
\big(\a+\frac{k+r}{k-s}(\b+r+k-1)\big).
\Eea 
Since $\overline{K(z,w)}^t=K(w,z)$, it follows that
$a_{mn}=\overline{a_{nm}}^t$ for $m,n\geq 0.$  Thus, by Lemma
\ref{relcoeff}, $(a_{01})_{r+1,r}=(r+1)!(\b)_{r+1}$, for $0\leq
r\leq n-1, (c_{k0}a_{00}^{-1})_{r,r+k}=\frac{(-1)^k}{k!}$, for
$0\leq r\leq n-k, 1\leq k\leq n+1$ and
\Bea
(c_{k0}a_{00}^{-1}a_{01})_{r,r+k-1}=(c_{k0}a_{00}^{-1})_{r,r+k}(a_{01})_{r+k,r+k-1}=
\frac{(-1)^k}{k!}(r+k)!(\b)_{r+k},
\Eea 
$0\leq r\leq n-k,1\leq k\leq n+1.$ Now, for $0\leq r\leq n-k, 2\leq k\leq n+1.$  Since
$c_{00}=a_{00}$, we clearly have 
\Bea
\lefteqn{(a_{00}^{1/2}\widetilde{\widetilde a}_{k1}a_{00}^{1/2})_{r,r+k-1}
=\bigg(\d_{s=1}^k c_{s0}a_{00}^{-1}a_{k-s,1}+a_{k1}\bigg)_{r,r+k-1}}\\
&&=\bigg(\d_{s=0}^{k-1}c_{s0}a_{00}^{-1}a_{k-s,1}+c_{k0}a_{00}^{-1}a_{01}\bigg)_{r,r+k-1}\\
&&=\d_{s=0}^{k-1}\tfrac{(-1)^s(k+r-1)!}{s!(k-s-1)!}(\b)_{r+k-1}\big(\a+\tfrac{k+r}{k-s}(\b+r+k-1)\big)
+\frac{(-1)^k(r+k)!}{k!}(\b)_{r+k}\\
&&=
\a(\b)_{r+k-1}\tfrac{(k+r-1)!}{(k-1)!}\d_{s=0}^{k-1}(-1)^s\tbinom{k-1}{s}+(\b)_{k+r}
\big(\d_{s=0}^{k-1}\tfrac{(-1)^s(k+r)!}{s!(k-s)!}+\tfrac{(-1)^k(k+r)!}{k!}\big)\\
&&=\tfrac{(k+r)!}{k!}(\bb S(z_2)^{\b+k+(k-l)}\b)_{k+r}\d_{s=0}^{k}(-1)^s\binom{k}{s}.
\Eea 
Therefore $(a_{00}^{1/2}\widetilde{\widetilde a}_{k1}a_{00}^{1/2})_{r,r+k-1}=0.$  Now,
$c_{00}=a_{00}$ and $(c_{k0}a_{00}^{-1}a_{01})_{n-k+1,n}=0$ for $2\leq k\leq n+1$.
Hence 
\Bea 
\lefteqn{(a_{00}^{1/2}\widetilde{\widetilde
a}_{k1}a_{00}^{1/2})_{n-k+1,n}
 = \bigg(\d_{s=1}^k c_{s0}a_{00}^{-1}a_{k-s,1}+a_{k1}\bigg)_{n-k+1,n}}\\
&&=\bigg(\d_{s=0}^{k-1}c_{s0}a_{00}^{-1}a_{k-s,1}\bigg)_{n-k+1,n}\\
&&=\d_{s=0}^{k-1}\tfrac{(-1)^s(k+(n-k+1)-1)!}{s!(k-s-1)!}(\b)_n
\bigg(\a+\tfrac{k+(n-k+1)}{k-s}(\b+n)\bigg)\\
&&=n!(\b)_{n}\bigg(\a\d_{s=0}^{k-1}\tfrac{(-1)^s}{s!(k-1-s)!}
+(n+1)(\b+n)\d_{s=0}^{k-1}\tfrac{(-1)^s}{s!(k-s)!}\bigg)\\
&&=n!(\b)_{n}\bigg(\tfrac{\a}{(k-1)!}\d_{s=0}^{k-1}(-1)^s\tbinom{k-1}{s}
+\tfrac{(n+1)(\b+n)}{k!}\d_{s=0}^k(-1)^s\tbinom{k}{s}-\tfrac{(-1)^k(n+1)(\b+n)}{k!}\bigg)\\
&&=0+0-n!(\b)_{n}\tfrac{(-1)^k(n+1)(\b+n)}{k!}\\
&&=\tfrac{(-1)^{k+1}(n+1)!(\b)_{n+1}}{k!},\mbox{~for~} 2\leq k\leq n+1.
\Eea 
Since $a_{00}$ is a diagonal matrix with positive diagonal entries, $\widetilde{\widetilde
a}_{k1}$ has the form as stated in the lemma, for $2\leq k\leq n+1,n\geq 1.$
\end{proof}
Here is a simple lemma from matrix theory which will be useful for us in the sequel.
\begin{lem}\label{mat}
Let $\{A_k\}_{k=0}^{n-1}$ are $(n+1)\times (n+1)$ matrices such
that $(A_k)_{kn}=\lambda_k\neq 0$ for $0\leq k\leq n-1,\,n\geq 1.$
If $AA_k=A_kA$ for some $(n+1)\times (n+1)$ matrix $A=\big ( \!\!
\big ( a_{ij}\big) \!\! \big )_{i,j=0}^n$ for $0\leq k\leq n-1,$ then $A$ 
is upper triangular with equal diagonal entries.
\end{lem}
\begin{proof}
$(AA_k)_{in}=a_{ik}(A_k)_{kn}=a_{ik}\l_k$ and
$(A_kA)_{kj}=(A_k)_{kn}a_{nj}=\l_ka_{nj}$ for $0\leq i,j\leq n,
0\leq k\leq n-1.$ Putting $i=k$ and $j=n,$   we have
$(AA_k)_{kn}=a_{kk}\l_k$ and $(A_kA)_{kn}=\l_ka_{nn}.$ By
hypothesis we have $a_{kk}\l_k=\l_ka_{nn}.$ As $\l_k\neq 0,$ this
implies that $a_{kk}=a_{nn}$ for $0\leq k\leq n-1,$ which is same
as saying that $A$ has equal diagonal entries. Now observe that
$(A_kA)_{ij}=0$ if $i\neq k$ for $0\leq j\leq n,$ which implies
that $(A_kA)_{in}=0$ if $i\neq k.$ By hypothesis this is same as
$(AA_k)_{in}=a_{ik}\l_{k}=0$ if $i\neq k$.  This implies
$a_{ik}=0$ if $i\neq k,0\leq i\leq n,0\leq k\leq n-1,$ which is a
stronger statement than saying $A$ is upper triangular.
\end{proof}
\begin{lem}
If an $(n+1)\times (n+1)$ matrix $A$ commutes with
$\widetilde{\widetilde a}_{k1}$ and $\widetilde{\widetilde
a}_{1k}$ for  $2\leq k\leq n+1,n\geq 1,$ then $A$ is a scalar.
\end{lem}
\begin{proof}
It follows from Lemma \ref{oncoeff} and Lemma \ref{mat} that if $A$ commutes with
$\widetilde{\widetilde a}_{k1}$ for $2\leq k\leq n+1,$ then $A$ is upper triangular
with equal diagonal entries. As the entries of $\widetilde{\widetilde a}_{k1}$ are
real, $(\widetilde{\widetilde a}_{1k})=(\widetilde{\widetilde a}_{k1})^t$.  If $A$
commutes with $\widetilde{\widetilde a}_{1k}$ for $2\leq k\leq n+1,$ then  by a
similar proof as in Lemma \ref{mat}, it follows that $A$ is lower triangular with
equal diagonal entries. So $A$ is both upper triangular and lower triangular with
equal diagonal entries, hence $A$ is a scalar.
\end{proof}
This sequence of Lemmas put together constitutes a proof of
Theorem 1.

For homogeneous operators in the class $\mathrm B_1(\D),$  we have a  
proof of reducibility that avoids the normalization of the kernel.  
This proof makes use of the fact that if such an operator is reducible then 
each of the direct summands must belong to the class $\mathrm B_1(\D).$ 
We give a precise formulation of this phenomenon along with a proof below. 
Let $K$ be a positive definite kernel on $\mathbb D^2$ and $\mathcal H$ 
be the corresponding Hilbert space.  Assume that the pair $(M_1,M_2)$ 
on $\mathcal H$ is in $\mathrm B_1(\D^2).$  The operator $M^*$ is the 
adjoint of the  multiplication operator on Hilbert space 
$J^{(2)}\mathcal H_{|{\rm res}~ \triangle}$ which consists of 
$\C^2$ - valued holomorphic function on $\D$ and possesses the reproducing kernel 
$J^{(2)}K(z,w)$. The operator $M^*$ is in $\mathrm B_2(\D)$ 
(cf. \cite[Proposition 3.6]{DMV}).
\begin{prop}\label{irr}
The operator $M^*$ on 
Hilbert space $J^{(2)}\mathcal H_{|{\rm res}~ \triangle}$ is irreducible.
\end{prop}

\begin{proof}
If possible, let $M^*$ be reducible, that is, $M^*=T_1\oplus T_2$ for some $T_1,
T_2\in\rm B_1(\D)$ by \cite[Proposition 1.18]{C-D}, which is same as saying by
\cite[Proposition 1.18]{C-D} that the associated bundle $E_{M^*}$ is reducible. A metric
on the associated bundle $E_{M^*}$ is given by $h(z)=J^{(2)}K(z,z)^t$.  
So, there exists a holomorphic change of frame
$\psi:\D\lo GL(2,\C)$ such that $\ov{\psi(z)}^th(z)\psi(z)=\left(%
\begin{array}{cc}
  h_1(z) & 0 \\
  0  & h_2(z) \\
\end{array}%
\right)$ for $z\in\D$, where $h_1$ and $h_2$ are metrics on the
associated line bundles $E_{T_1}$
and $E_{T_2}$ respectively. So ${\psi(z)}^{-1}\mathcal K_h(z)\psi(z)=\left(%
\begin{array}{cc}
  \mathcal K_{h_1}(z) & 0 \\
  0 & \mathcal K_{h_2}(z)\\
\end{array}%
\right)$, where $\mathcal K_h(z)=\bar\partial(h^{-1}\partial h)(z)$ is the curvature of the
bundle $E_{M^*}$ with respect to the metric $h$ and  
 $\mathcal K_{h_i}(z)$ are the curvatures of the bundles $E_{T_i}$ for $i=1,2$ as in 
\cite[pp. 211]{C-D}. 
A direct computation shows that $\K_h(z)=\left(%
\begin{array}{cc}
  \a & -2\b(\b+1)(1-|z|^2)^{-1}\bar z \\
  0 & \a+2\b+2 \\
\end{array}%
\right)(1-|z|^2)^{-2}.$  Thus the matrix $\psi(z)$  diagonalizes
$\K_h(z)$ for $z\in\D$.  It follows that $\psi(z)$ is determined, that is, 
the columns of $\psi(z)$ are eigenvectors of $\K_h(z)$ for
$z\in\D.$  These are uniquely determined upto multiplication by
non-vanishing scalar valued functions $f_1$ and $f_2$ on $\D$. Now
one set of eigenvectors of $\K_h(z)$ is given
by $\{\left(%
\begin{array}{c}
  1 \\
  0 \\
\end{array}%
\right),\left(%
\begin{array}{c}
  -\b\bar z \\
  1-|z|^2 \\
\end{array}%
\right)\}$ and it is clear that there does not exist any non-vanishing scalar valued
function $f_2$ on $\D$ such that  $f_2(z)\left(%
\begin{array}{c}
  -\b\bar z \\
  1-|z|^2 \\
\end{array}%
\right)$  is an eigenvector for $\K_h(z)$ whose entries are
holomorphic functions on $\D$. Hence there does not exist any
holomorphic change of frame $\psi:\D\lo
GL(2,\C)$ such that $\ov{g}^thg=\left(%
\begin{array}{cc}
  h_1 & 0 \\
  0  & h_2 \\
\end{array}%
\right)$ on $\D.$  Hence $M^*$ is irreducible.
\end{proof}
\begin{thm}
The operators $T=M^{(\a,\b)}$ and $\wi T:=M^{(\wi\a,\wi \b)}$ are unitarily
equivalent if and only if $\a=\wi\a$ and $\b=\wi\b.$
\end{thm}
One of the implications is trivial.  To prove the other implication, recall that 
\cite[Proposition 3.6]{DMV}   $T,\wi T\in\mathrm B_{n+1}(\D)$.  It follows from
\cite{C-D} that if $T ,\wi T\in\mathrm B_{n+1}(\D)$ are unitarily equivalent then
the curvatures $\K_T,\K_{\wi T}$ of the associated bundles $E_T$ and $E_{\wi T}$
respectively, are unitarily equivalent as matrix-valued real-analytic functions on
$\D$. In particular, this implies that $\K_T(0)$ and $\K_{\wi T}(0)$ are \u.
Therefore, we compute $\K_T(0)$ and $\K_{\wi T}(0).$  Let $\K_T(h)$ denote 
the curvature of the bundle $E_T$ with respect to the metric 
$h(z):=\wi{\wi K}(z,z)^t$. 
 
\begin{lem}\label{cur}
The curvature $\K_T(h)(0)$ at $0$ of the bundle $E_T$ equals the coefficient 
of $z \bar{z}$ in the normalized kernel $\wi{\wi K}$, that is,  $\K_T(h)(0)=\wi{\wi a}_{11}^t.$ 
\end{lem}
\begin{proof}
The curvature of the bundle $E_T$  with respect to the metric $h(z):=\wi{\wi K}(z,z)^t$  
is  $\K_T(h)=\bar\partial(h^{-1}\partial h)$. If $h(z)=\sum_{m,n\geq
0}h_{mn}z^m{\bar z}^n$, then $h_{mn}={\wi{\wi a}_{mn}^t}$ for
$m,n\geq 0.$ So, $h_{00}=I$ and $h_{m0}=h_{0n}=0$ for $m,n\geq 1.$
Hence $\K_T(h)(0)=\bar\partial h^{-1}(0)\partial
h(0)+h^{-1}(0)\bar\partial\partial h(0)=(\bar\partial
h^{-1}(0))h_{10}+h_{00}^{-1}h_{11}=h_{11}=\wi{\wi a}_{11}^t.$
\end{proof}
\begin{lem}
$\big(\K_T(0)\big)_{ii}=\a,$ for $i=0,\ldots, n-1$ and
$\big(\K_T(0)\big)_{nn}=\a+(n+1)(\b+n)$ for $n\geq 1.$
\end{lem}
\begin{proof}
From Lemma \ref{cur} and Lemma \ref{coeff} we know that $\K_T(0)=\wi{\wi a}_{11}^t
=\big(\wi a_{11}+\widetilde b_{10}\wi a_{01}\big)^t.$ Thus $\K_T(0)$ is the
transpose of $a_{00}^{-1/2}(a_{11}+c_{10}a_{00}^{-1}a_{01})a_{00}^{-1/2}$ by Lemma
\ref{impcoeff}. Now, by Lemma \ref{relcoeff} and Lemma \ref{impcoeff},
$\big(c_{10}\big)_{r,r+1}=-(r+1)!(\b)_{r+1}$ for $0\leq r\leq n-1,$
$\big(a_{00}\big)_{rr}=r!(\b)_r,
\big(a_{11}\big)_{rr}=r!(\b)_r\big(\a+(r+1)(\b+r)\big)$ for $0\leq r \leq n$ and
$\big(a_{01}\big)_{r+1,r}=(r+1)!(\b)_{r+1}$ for $0\leq r\leq n-1.$ Therefore, 
$\big(c_{10}a_{00}^{-1}a_{01}\big)_{rr}=-(r+1)!(\b)_{r+1}$ for $0\leq r\leq n-1.$
Also, $\big(a_{11}+c_{10}a_{00}^{-1}a_{01}\big)_{rr}=\a r!(\b)_{r+1}$ for $0\leq r\leq
n-1,$ and $\big(a_{11}+c_{10}a_{00}^{-1}a_{01}\big)_{nn}=n!(\b)_n(\a+(n+1)(\b+n))$. Finally, 
$\K_T(h)(0)=\wi{\wi a}_{11}^t=\wi{\wi a}_{11},$ as $\wi{\wi a}_{11}$ is a diagonal
matrix with real entries. In fact, $\big(\K_T(0)\big)_{ii}=\a,$ for $i=0,\ldots, n-1$
and $\big(\K_T(0)\big)_{nn}=\a+(n+1)(\b+n).$
\end{proof}
We now see that $M$ and $\wi M$ are \u\  implies  that  $\a=\wi\a$ and
$\a+(n+1)(\b+n)=\wi\a+(n+1)(\wi\b+n)$, that is, $\a=\wi\a$ and $\b=\wi\b.$
\section{Homogeneity of the operator $M^{(\a,\b)}$}

\begin{thm}\label{homo}
The multiplication operator $M:=M^{(\a,\b)}$ on $J^{(n+1)}\mathcal
H$ is homogeneous.
\end{thm}
This theorem is a particular case of the Lemma \ref{newpair}.  A proof first
appeared in \cite[Theorem 5.2.]{BMIAS01}. We give an alternative proof of this
Theorem by showing that that the kernel is quasi-invariant, that is,
$$K(z,w)=J_{\p}(z)K\big(\p(z),\p(w)\big)\overline{J_{\p}(w)}^t$$ 
for some cocycle 
$$J:\mbox{\rm M\"{o}b}\times \D\lo \C^{(n+1)\times (n+1)},\, \v \in \mbox{\rm M\"{o}b},\, z,\, 
w \in\D.$$ 
First we prove that $K(z,z)=J_{\p}(z)K\big(\p(z),\p(z)\big)\overline{J_{\p}(z)}^t$ 
and then polarize to obtain the final result.  We begin with a series of lemmas.
\begin{lem} \label{homog}
Suppose that $J: \mbox{\rm M\"{o}b}\times \D\lo \C^{(n+1)\times (n+1)}$ is a
cocycle. Then the following are equivalent
\begin{enumerate}
\item $K(z,z)=J_{\p}(z)K\big(\p(z),\p(z)\big)\overline{J_{\p}(z)}^t$ for all $\v\in$
{M\"{o}b} and $z\in\D;$ 
\item $K(0,0)=J_{\p}(0)K\big(\p(0),\p(0)\big)\overline{J_{\p}(0)}^t$ for all $\v\in$
{M\"{o}b}.
\end{enumerate}
\end{lem}
\begin{proof}
One of the implications is trivial. To prove the other implication, note that 
\Bea
J_{\p_1}(0)K\big(\p_1(0),\p_1(0)\big)\overline{J_{\p_1}(0)}^t&=&K(0,0)\\
&=&J_{\p_2}(0)K\big(\p_2(0),\p_2(0)\big)\overline{J_{\p_2}(0)}^t
\Eea
for any $\v_1,\v_2\in${M\"{o}b} and $z\in\D.$  Now pick
$\psi\in$M\"{o}b such that $\s(0)=z$ and taking $\v_1=\psi,
\v_2=\psi\v$ in the previous identity we see that
\begin{eqnarray*}
\lefteqn{J_{\s}(0)K\big(\s(0),\s(0)\big)\overline{J_{\s}(0)}^t}\\
&=& J_{\p\s}(0)K\big(\p\s(0),\p\s(0)\big)\overline{J_{\p\s}(0)}^t\\
&=& J_{\s}(0)J_{\p}(\s(0))K\big(\p(z),\p(z)\big)\overline{J_{\p}(\s(0))}^t
\overline{J_{\s}(0)}^t
\end{eqnarray*}
for $\v\in\mbox{M\"{o}b}, z\in\D.$  Since $J_{\s}(0)$ is
invertible, it follows from the equality of first and third quantities that
$$K\big(\s(0),\s(0)\big)=J_{\p}(\s(0))K\big(\p\s(0),\p\s(0)\big)\overline{J_{\p}(\s(0))}^t.$$
This is the same as $K(z,z)=J_{\p}(z)K\big(\p(z),\p(z)\big)\overline{J_{\p}(z)}^t$ by
the choice of $\psi.$  The proof of this lemma is therefore complete.
\end{proof}
Let $\j_{\p}(z)=(J_{\p}(z)^t)^{-1},\v\in$M\"{o}b, $z\in\D$, where $X^t$ denotes the
 transpose of the matrix $X$. Clearly, $(J_{\p}(z)^t)^{-1}$ satisfies the cocycle
property if and only if $\j_{\p}(z)$ does and they uniquely determine each other. It
is easy to see that the condition $$K(0,0)=J_{\p}(0)K\big(\p(0),\p(0)\big)\overline{J_{\p}(0)}^t$$
is equivalent to
\begin{equation} \label{cocyclerule}
h\big(\p(0)\big)=\overline{\j_{\p}(0)}^th(0)\j_{\p}(0),
\end{equation}
where $h(z)$  is the transpose of  $K(z,z)$ as before.  It will be useful to define the two functions   
\begin{enumerate}
\item[(i)] $c:$ M\"{o}b$\times \D\lo\C$, $c(\p,z)=(\p)^\prime(z)$; 
\item[(ii)] $p:$ M\"{o}b$\times \D\lo\C$, $p(\p,z)=\frac{\ov{ta}}{1+\ov{ta}z}$ 
\end{enumerate}
for $\v_{t,a}\in  \mbox{\rm M\"{o}b},\,t \in\mathbb T,\,a \in  \D.$  We point out that the 
function $c$ is the well-known cocycle for the group  M\"{o}b.  
\begin{lem} \label{imp}
With notation as above, we have 
\begin{enumerate}
\item[(a)] $\p_{t,a}=\v_{\bar t,-ta}$
\item[(b)]$\v_{s,b} \v_{t,a}=\v_{\frac{s(t+\bar ab)}{1+ta\bar b},\frac{a+\bar
tb}{1+\ov{ta}b}}$ 
\item[(c)] $c(\p,\s(z))c(\s(z))=c(\p \s,z)$ for $\v,\psi\in$\rm{M\"{o}b}, $z\in\D$ 
\item[(d)]  $p(\p,\s(z))c(\s,z)+p(\s,z)=p(\p\s,z)$ for $\v,\psi\in$\rm{M\"{o}b}, $z\in\D.$
\end{enumerate}
\end{lem}
\begin{proof}
The proof of (a) is a mere  verification.
We note that $$ \v_{s,b}(\v_{t,a}(z))=s\frac{t\frac{z-a}{1-\bar az}-b}{1-\bar
bt\frac{z-a}{1-\bar az}}=s\frac{tz-ta-b+\bar abz}{1-\bar az-t\bar bz+ta\bar
b}=\frac{s(t+\bar ab)}{1+ta\bar b}\frac{z-\frac{ta+b}{t+\bar ab}}{1-\frac{\bar
a+t\bar b}{1+ta\bar b}z},$$ which is (b).
The chain rule gives (c).
To prove (d), we first note that  for $\v=\v_{t,a}$ and $\psi=\v_{s,b}$,  if 
$\psi^{-1}\v^{-1} = \v_{t^\prime,a^\prime}$ for some $(t^\prime,a^\prime)\in\mathbb T\times \D$
then 
$$\ov{t^\prime a^\prime}=\frac{\bar s(\bar t+a\bar b)}{1+\ov{ta}b}\frac{\bar a+t\bar b}{1+ta\bar b} =
\frac{\bar s(\bar b+\ov{ta})}{1+\ov{ta}b}.$$
It is now easy to verify that  
\Bea
p(\p,\s(z))c(\s,z)+p(\s,z)
&=&\frac{\ov{ta}}{1+\ov{ta}\s_{s,b}(z)}\frac{\bar
s(1-|b|^2)}{(1+\ov{sb}z)^2}+\frac{\ov{sb}}{1+\ov{sb}z}\\
&=&\frac{\bar s(\bar b+\ov{ta})}{1+\ov{ta}b+\bar s(\bar
b+\ov{ta})z}\\
&=&\frac{\bar s\frac{\bar b+\ov{ta}}{1+\ov{ta}b}}{1+\frac{\bar
s(\bar b+\ov{ta})}{1+\ov{ta}b}z}\\
&=&p(\p\s,z). 
\Eea
\end{proof}
\noindent 
Let
\begin{equation}  \label{calJ}
\big(\j_{\p}(z)\big)_{ij}=c(\p,z)^{-\frac{\a+\b}{2}-n}\frac{(\b)_j}{(\b)_i}\binom{j}{i}
c(\p,z)^{n-j}p(\p,z)^{j-i}
\end{equation} 
for $0\leq i\leq j\leq n.$ 
\begin{lem} \label{L18}
$J_{\p}(z)$ defines a cocycle for the group M\"{o}b.
\end{lem}
\begin{proof}
To say that $J_{\p}(z)$ satisfies the cocycle property  is the same as 
saying $\j_{\p}(z)$ satisfies the cocycle property, which is what we will 
verify.  Thus we want to show that
$\big(\j_{\s}(z)\j_{\p}(\s(z))\big)_{ij}=\big(\j_{\p\s}(z)\big)_{ij}$
for $0\leq i,j\leq n$. We note that $\j_{\p}(z)$ is upper
triangular, as the product of two upper triangular matrices is
again upper triangular, it suffices  to prove this equality for $0\leq
i\leq j\leq n.$  Clearly, we have 
\Bea 
\lefteqn{\big(\j_{\s}(z)\j_{\p}(\s(z))\big)_{ij} =\displaystyle\sum_{k=i}^j\big(\j_{\s}(z)\big)_{ik}
\big(\j_{\p}(\s(z))\big)_{kj} }\\
&& =c(\s,z)^{-\frac{\a+\b}{2}-n}c(\p,\s(z))^{-\frac{\a+\b}{2}-n}
\d_{k=i}^j \Big ( \frac{(\b)_k}{(\b)_{i}}\binom{k}{i}c(\s,z)^{n-k}\\
&&p(\s,z)^{k-i}\frac{(\b)_j}{(\b)_k}\binom{j}{k}c(\p,\s(z))^{n-j}p(\p,\s(z))^{j-k} \Big)\\
&&=c(\p\s,z)^{-\frac{\a+\b}{2}-n}\frac{(\b)_j}{(\b)_i}c(\s,z)^{n-j}c(\p,\s(z))^{n-j}\\
&&\hspace*{36pt}\d_{k=i}^j\frac{j!}{i!(k-i)!(j-k)!}c(\s,z)^{j-k}p(\p,\s(z))^{j-k}p(\s,z)^{k-i}\\
&&=c(\p\s,z)^{-\frac{\a+\b}{2}-n}\frac{(\b)_j}{(\b)_i}\binom{j}{i}c(\p\s,z)^{n-j}\\
&&\hspace*{36pt}\d_{k=i}^j\binom{j-i}{k-i}c(\s,z)^{j-k}p(\p,\s(z))^{j-k}p(\s,z)^{k-i}\\
&&=c(\p\s,z)^{-\frac{\a+\b}{2}-n}\frac{(\b)_j}{(\b)_i}\binom{j}{i}c(\p\s,z)^{n-j}\\
&&\phantom{Gadadhar}\d_{k=0}^{j-i}\binom{j-i}{k}c(\s,z)^{(j-i)-k}
p(\p,\s(z))^{(j-i)-k}p(\s,z)^{k}\\
&&=c(\p\s,z)^{-\frac{\a+\b}{2}-n}\frac{(\b)_j}{(\b)_i}\binom{j}{i}c(\p\s,z)^{n-j}\\
&&\phantom{GadadharGadadhar}
\bigg(c(\s,z)p(\p,\s(z))+p(\s,z)\bigg)^{j-i}\\
&&=c(\p\s,z)^{-\frac{\a+\b}{2}-n}\frac{(\b)_j}{(\b)_i}\binom{j}{i}c(\p\s,z)^{n-j}
p(\p\s,z)^{j-i}\\
&&=\big(\j_{\p\s}(z)\big)_{ij},\Eea
for  $0\leq i\leq j\leq n$.  The penultimate equality follows from
Lemma \ref{imp}.
\end{proof}

We need the following beautiful identity to prove \mbox{(\ref{cocyclerule})}. We
provide two proofs, the first one is due to C. Verughese and the second is due to B.
Bagchi.

\begin{lem}\label{iden}
For nonnegative integers $j\geq i$ and $0\leq k\leq i$, we have 
 \begin{multline*}
 \sum_{l=0}^{i-k}(-1)^l(l+k)!\binom{i}{l+k}\binom{j}{l+k}\binom{l+k}{l}(a+j)_{i-l-k}
 = k!\binom{i}{k}\binom{j}{k}(a+k)_{i-k},
\end{multline*}
for all $a\in\C$.
\end{lem}
\begin{proof} 
Here is the first proof due to C. Verughese:  
For any integer $i\geq 1$ and  $a\in \mathbb C \setminus \mathbb Z,$ we have
\Bea
\lefteqn{\sum_{l=0}^{i-k}(-1)^l
(l+k)!\binom{i}{l+k}\binom{j}{l+k}\binom{l+k}{l}(a+j)_{i-l-k}}\\
&=&
\frac{i!j!}{k!\Gamma(a+j)}\sum_{l=0}^{i-k}\frac{(-1)^l}{l!(i-k-l)!}
\frac{\Gamma(a+j+i-l-k)}{\Gamma(j-l-k+1)}\\
&=& \frac{i!j!}{k!(i-k)!\Gamma(a+j)\Gamma(1-a-i)}
\sum_{l=0}^{i-k}(-1)^l\tbinom{i-k}{l}B(a+j+i-k-l, 1-a-i)\\
&=& \frac{i!j!}{k!(i-k)!\Gamma(a+j)\Gamma(1-a-i)}
\sum_{l=0}^{i-k}(-1)^l\tbinom{i-k}{l}\int_{0}^1
t^{a+j+i-k-l-1}(1-t)^{-a-i}dt\\
&=& \frac{i!j!}{k!(i-k)!\Gamma(a+j)\Gamma(1-a-i)}
\int_{0}^1\sum_{l=0}^{i-k}(-1)^l\tbinom{i-k}{l}
t^{a+j+i-k-l-1}(1-t)^{-a-i}dt\\
&=& \frac{i!j!}{k!(i-k)!\Gamma(a+j)\Gamma(1-a-i)}
\int_{0}^1(1-t)^{-a-i}t^{a+j-1}\big(\sum_{l=0}^{i-k}(-1)^l\tbinom{i-k}{l}t^{i-k-l}\big)dt\\
&=& \frac{i!j!}{k!(i-k)!\Gamma(a+j)\Gamma(1-a-i)}
\int_{0}^1(1-t)^{-a-i}t^{a+j-1}(t-1)^{i-k}dt\\
&=&\frac{(-1)^{i-k}i!j!}{k!(i-k)!\Gamma(a+j)\Gamma(1-a-i)}
\int_{0}^1(1-t)^{-a-i}t^{a+j-1}(1-t)^{i-k}dt\\
&=& \frac{(-1)^{i-k}i!j!}{k!(i-k)!\Gamma(a+j)\Gamma(1-a-i)}B(a+j,
1-a-k)\\
&=& \frac{(-1)^{i-k}i!j!}{k!(i-k)!\Gamma(a+j)\Gamma(1-a-i)}
\frac{\Gamma(a+j)\Gamma(1-a-k)}{\Gamma(1+j-k)}\\
&=& \frac{(-1)^{i-k}i!j!}{k!(i-k)!\Gamma(1-a-i)}
\frac{\Gamma(1-a-k)}{(j-k)!}\\
&=&
k!\binom{i}{k}\binom{j}{k}\frac{\Gamma(1-a)}{(-1)^i\Gamma(1-a-i)}
\frac{(-1)^k\Gamma(1-a-k)}{\Gamma(1-a)}\\
&=& k!\binom{i}{k}\binom{j}{k} \frac{\Gamma(1-a)}{\cos
i\pi\Gamma(1-a-i)} \frac{\cos
k\pi\Gamma(1-a-k)}{\Gamma(1-a)}\\
&=& k!\binom{i}{k}\binom{j}{k} \frac{\Gamma(1-a)\Gamma(a+i)}{\cos
i\pi\Gamma(1-a-i)\Gamma(a+i)} \frac{\cos
k\pi\Gamma(1-a-k)\Gamma(a+k)}{\Gamma(a+k)\Gamma(1-a)}\\
&=& k!\binom{i}{k}\binom{j}{k}
\frac{\Gamma(1-a)\Gamma(a+i)\sin(a+i)\pi}{\pi\cos (i\pi)}
\frac{\pi\cos
k\pi}{\sin (a+k)\pi\Gamma(a+k)\Gamma(1-a)}\\
&=&
k!\binom{i}{k}\binom{j}{k}\frac{\Gamma(a+i)}{\Gamma(a+k)}.\Eea

Since we have an equality involving a polynomial of degree $i-k$ for all $a$
in $\mathbb C \setminus \mathbb Z$, it follows that the equality holds
for all $a\in \mathbb C$.

Here is another proof due to B. Bagchi: 
Since $\binom{-x}{n}=\frac{-x(-x-1)\cdots(-x-n+1)}{n!}
=(-1)^n\binom{x+n-1}{n}$ and $(x)_n=x(x+1)\cdots(x+n-1)=n!\binom{x+n-1}{n},$ 
it follows that  
\Bea
\lefteqn{\sum_{l=0}^{i-k}(-1)^l
(l+k)!\binom{i}{l+k}\binom{j}{l+k}\binom{l+k}{l}(a+j)_{i-l-k}}\\
&=&
\frac{i!j!}{k!}\sum_{l=0}^{i-k}\frac{(-1)^l}{l!(i-k-l)!(j-k-l)!}
(i-k-l)!\binom{a+j+i-k-l-1}{i-k-l}\\
&=&
\frac{i!j!}{k!(j-k)!}\sum_{l=0}^{i-k}\frac{(-1)^l(j-k)!}{l!(j-k-l)!}
(-1)^{i-k-l}\binom{-a-j}{i-k-l}\\
&=&
i!\binom{j}{k}(-1)^{i-k}\sum_{l=0}^{i-k}\binom{j-k}{l}\binom{-a-j}{i-k-l}\\
&=& i!\binom{j}{k}(-1)^{i-k}\binom{-a-k}{i-k}\\
&=& i!\binom{j}{k}(-1)^{i-k}(-1)^{i-k}\binom{a+i-1}{i-k}\\
&=& k!\binom{i}{k}\binom{j}{k}(a+k)_{i-k},
\Eea 
where the equality after the last summation symbol follows from Vandermonde's
identity.
\end{proof}

\begin{lem} For $\v\in${M\"{o}b} and $\j_{\p}(z)$ as in  (\ref{calJ}),
$$h\big(\p(0)\big)=\overline{\j_{\p}(0)}^th(0)\j_{\p}(0).$$

\end{lem}
\begin{proof}
It is enough to show that 
\Bea
h\big(\p(0)\big)_{ij}=\big(\overline{\j_{\p}(0)}^t
h(0)\j_{\p}(0)\big)_{ij} , \mbox{~for~} 0\leq i\leq j\leq n.
\Eea 
Let $\v=\v_{t,z},\, t\in\mathbb T,$ and $z\in\D.$  Since
$\big(h\big(\p(0)\big)\big)_{ij}=\big(h(z)\big)_{ij}$, it follows that
\Bea
\lefteqn{
\big(h\big(\p(0)\big)\big)_{ij}
={\bar\partial}_2^i\partial_2^j
\big(\bb S(z_1)^\a\bb S(z_2)^\b\big)|_{\triangle}}\\
&&=(\b)_j\bb S(z_1)^\a{\bar\partial}_2^i\big(\bb S(z_2)^{\b+j}\bar
z_2^j\big)|_{\triangle}\\
&&=
(\b)_j\bb S(z_1)^\a\d_{r=0}^i\binom{i}{r}{\bar\partial}_2^{(i-r)}\big(\bb S(z_2)^{\b+j}\big)
{\bar\partial}_2^r(\bar z_2^j)|_{\triangle}\\
&&=(\b)_j\bb S(z_1)^\a\d_{r=0}^i\binom{i}{r}(\b+j)_{i-r}\bb S(z_2)^{\b+j+(i-r)}
z_{2}^{i-r}r!\binom{j}{r}\bar z_2^{j-r}|_{\triangle}\\
&&=(\b)_j\bb S(z)^{\a+\b+i+j}{\bar
z}^{j-i}\d_{r=0}^ir!\binom{i}{r}\binom{j}{r}(\b+j)_{i-r}
\bb S(z)^{-r} |z|^{2(i-r)},
\Eea 
for $i\leq j.$ Clearly, $\big(\j_{\p}(0)\big)_{ij}=c(\p,0)^{-\frac{\a+\b}{2}-n}\frac{(\b)_j}{(\b)_i}\binom{j}{i}
c(\p,0)^{n-j}p(\p,0)^{j-i}$ and $h(0)_{ii}=i!(\b)_i,$ $0\leq i\leq j\leq
n.$  We have 
\Bea 
\big(\overline{\j_{\p}(0)}^th(0)\j_{\p}(0)\big)_{ij} &=&
\d_{k=0}^j\big(\overline{\j_{\p}(0)}^th(0)\big)_{ik}\big(\j_{\p}(0)\big)_{kj}\\
&=&\d_{k=0}^i\d_{k=0}^j\big(\overline{\j_{\p}(0)}^t\big)_{ik}
\big(h(0)\big)_{kk}\big(\j_{\p}(0)\big)_{kj}\\
&=&\d_{k=0}^{\mbox{min}(i,j)}\big(\overline{\j_{\p}(0)}^t\big)_{ik}
\big(h(0)\big)_{kk}\big(\j_{\p}(0)\big)_{kj}.\Eea Now, for
$0\leq i\leq j\leq n,$ 
\begin{multline*}
\d_{k=0}^{\mbox{min}(i,j)}\big(\overline{\j_{\p}(0)}^t\big)_{ik}
\big(h(0)\big)_{kk}\big(\j_{\p}(0)\big)_{kj}=|c(\p,0)|^{-\a-\b-2n}\\
\d_{k=0}^i\Big ( \frac{(\b)_i}{(\b)_{k}}
\binom{i}{k}\ov{c(\p,0)}^{n-i}\ov{p(\p,0)}^{i-k}
k!(\b)_k\frac{(\b)_j}{(\b)_{k}}\\
\phantom{GadadharGadadharGadadhar}
\binom{j}{k}{c(\p,0)}^{n-j}{p(\p,0)}^{j-k} \Big )\\
=\bb S(z)^{\a+\b+2n}\d_{k=0}^i\frac{k!(\b)_{i}(\b)_{j}}{(\b)_{k}}\binom{i}{k}\binom{j}{k}\\
\phantom{GadadharGadadharGadadharGadadhar}
\big(t\bb S(z)\big)^{-n+i}(tz)^{i-k}\big(\ov t\bb S(z)\big)^{-n+j}(\ov{tz})^{j-k}\\
=(\b)_j\bb S(z)^{\a+\b+i+j}\d_{k=0}^ik!\binom{i}{k}\binom{j}{k}
\frac{(\b)_{i}}{(\b)_{k}}z^{i-k}\bar z^{j-k}\\
=(\b)_j\bb S(z)^{\a+\b+i+j}\bar z^{j-i}\d_{k=0}^ik!\binom{i}{k}\binom{j}{k}
\frac{(\b)_{i}}{(\b)_{k}}|z|^{2(i-k)}.
\end{multline*}
Clearly, to prove the desired equality we have to show that 
\bea 
\d_{r=0}^ir!\binom{i}{r}\binom{j}{r}(\b+j)_{i-r}
\bb S(z)^{-r}|z|^{2(i-r)}=
\d_{k=0}^ik!\binom{i}{k}\binom{j}{k}
\frac{(\b)_{i}}{(\b)_{k}}|z|^{2(i-k)} \label{equality}
\eea 
$\mbox{~for~}0\leq i\leq j\leq n.$ But 
\Bea 
\lefteqn{\d_{r=0}^ir!\binom{i}{r}\binom{j}{r}(\b+j)_{i-r}
(1-|z|^2)^r|z|^{2(i-r)}}\\
&&=\d_{r=o}^ir!\binom{i}{r}\binom{j}{r}(\b+j)_{i-r}
\d_{l=0}^r(-1)^l\binom{r}{l}|z|^{2l}|z|^{2(i-r)}\\
&&=\d_{l=0}^i\d_{r=l}^i
(-1)^lr!\binom{i}{r}\binom{j}{r}\binom{r}{l}(\b+j)_{i-r}|z|^{2(i-(r-l))}\\
&&=\d_{l=0}^i\d_{r=0}^{i-l}
(-1)^l(r+l)!\binom{i}{r+l}\binom{j}{r+l}\binom{r+l}{l}(\b+j)_{i-r-l}|z|^{2(i-r)}.
\Eea
For $0\leq k\leq i-l,$ the coefficient of $|z|^{2(i-k)}$ in the
left hand side of (\ref{equality}) is
$$\d_{l=0}^i(-1)^l(k+l)!\binom{i}{k+l}\binom{j}{k+l}\binom{k+l}{l}(\b+j)_{i-k-l},$$
which is the same as
$$\d_{l=0}^{i-k}(-1)^l(k+l)!\binom{i}{k+l}\binom{j}{k+l}\binom{k+l}{l}(\b+j)_{i-k-l},$$
for $0\leq l\leq i-k\leq i.$  So, to complete the proof we have to
show that
$$\d_{l=0}^{i-k}(-1)^l(k+l)!\binom{i}{k+l}\binom{j}{k+l}\binom{k+l}{l}(\b+j)_{i-k-l}
=k!\binom{i}{k}\binom{j}{k} \frac{(\b)_{i}}{(\b)_{k}},$$ for $0\leq
k\leq i,i\leq j.$ But this follows from Lemma \ref{iden}.
\end{proof}

\section{The case of the tri-disc $\mathbb D^3$}
We discuss the jet construction for $\D^3$. Let $K:\D^3\times\D^3\lo\C$ be a
reproducing kernel.  Following the jet construction of \cite{DMV}, we define  
$$J^{(1,1)}K(z,w)=\left(%
\begin{array}{ccc}
  K(z,w) & \partial_2K(z,w) &  \partial_3K(z,w) \\
  \bar\partial_2K(z,w) &  \partial_2\bar\partial_2K(z,w) &  \bar\partial_2\partial_3K(z,w) \\
  \bar\partial_3K(z,w)  & \partial_2\bar\partial_3K(z,w) & \bar\partial_3\partial_3K(z,w) \\
\end{array}\right),\,\, z,\,w\in \D^3.$$%
As before, to
retain the usual meaning of $\partial$ and $\bar\partial$ we replace
$J^{(1,1)}K(z,w)$ by its transpose. For simplicity of notation, we let $G(z,w):=J^{(1,1)}K(z,w)^t$.
In this notation, choosing the kernel function $K$ on $\D^3$ to be  $$K(z,w)=(1-z_1\bar{w}_1)^{-\alpha}(1-z_2\bar{w}_2)^{-\beta}(1-z_3\bar{w}_3)^{-\gamma},$$ 
we have 
$$G(z,w)=\left(%
\begin{array}{ccc}
  (1-z\bar w)^2 & \b z(1-z\bar w) & \gamma z(1-z\bar w) \\
 \b\bar w(1-z\bar w)& \b(1+\b z\bar w) & \b\gamma z\bar w \\
  \gamma \bar w(1-z\bar w) & \b\gamma z\bar w & \gamma(1+\gamma z \bar w) \\
\end{array}%
\right)(1-z\bar w)^{-\a-\b-\gamma-2},$$ for $z,w\in\D,~~~\a,\b,\g>0.$
\begin{thm}
The adjoint of the multiplication operator $M^*$ on the
Hilbert space  of $\C^3$ valued holomorphic functions on $\D^3$ with
reproducing kernel $G$ is in $\mathrm B_3(\D)$. It is homogeneous and reducible.   
Moreover, $M^*$ is \u\ to
$M_1^*\oplus M_2^*$ for a pair of irreducible homogeneous operators 
$M_1^*$ and $M_2^*$ from $\mathrm B_1(\D)$. 
\end{thm}
\begin{proof}
Although homogeneity of $M^*$ follows from \cite[Theorem 5.2.]{BMIAS01}, we
give an independent proof using the ideas we have developed in this note.
Let $$\wi{\wi G}(z,w)=G(0,0)^{1/2}G(z,0)^{-1}G(z,w)G(0,w)^{-1}G(0,0)^{1/2}.$$ 
Evidently, $\wi {\wi G}(z,0)=I,$ that is, $\wi{\wi G}$ is a normalized kernel. The form of
$\wi {\wi G}(z,w)$ for $z,w\in\D$ is $(1-z\bar w)^{-\a-\b-\g-2}$ times 
$$
\left ( \begin{smallmatrix}  
  (1-z\bar w)^2 -(\b+\g)(1-z\bar w)z\bar w\\
  +(\b+\g)(1+\b+\g)z^2{\bar w}^2 & -\sqrt\b(1+\b+\g)z^2\bar w & -\sqrt\g(1+\b+\g)z^2\bar w\\
 &\\
 -\sqrt\b(1+\b+\g)z{\bar w}^2& 1+\b z\bar w & \sqrt{\b\g}z\bar w \\
 &\\
    -\sqrt\g(1+\b+\g)z{\bar w}^2 & \sqrt{\b\g}z\bar w  & 1+\g z\bar w\\
\end{smallmatrix}
\right).$$
Let $U= \frac{1}{\sqrt{\b+\g}} \left(%
  \begin{smallmatrix}
  1 & 0 & 0 \\
  0 & \sqrt{\b} & \sqrt{\g} \\
  0 & -\sqrt{\g} & \sqrt{\b} \\
\end{smallmatrix}
\right)$  
which is unitary on $\C^3$.  By a direct computation, we see that the
equivalent normalized kernel $U\wi {\wi G}(z,w)\ov U^t$ is equal to the direct sum
$G_1(z,w) \oplus G_2(z,w)$, where
$G_2(z,w)= (1-z\bar w)^{-\a-\b-\g-2}$ and $$G_1(z,w)=\left(%
  \begin{smallmatrix}
  (1-z\bar w)^2 -(\b+\g)(1-z\bar w)z\bar w\\
  +(\b+\g)(1+\b+\g)z^2{\bar w}^2 & -\sqrt{\b+\g}(1+\b+\g)z^2\bar w \\
  & \\
  -\sqrt{\b+\g}(1+\b+\g)z{\bar w}^2  & 1+(\b+\g)z\bar w  \\
  \end{smallmatrix}
\right)(1-z\bar w)^{-\a-\b-\g-2}.$$ It follows that $M^*$ is \u\ to a
reducible operator  by an application of Lemma \ref{kur}, that is, $M^*$ is
reducible. If we replace $\b$ by $\b+\g$ in Theorem
\ref{tirr} take $n=1,$ then $$K(z,w)=\left(%
\begin{array}{cc}
  (1-z\bar w)^2 & (\b+\g)z(1-z\bar w) \\
  (\b+\g)\bar w(1-z\bar w) & (\b+\g)(1+(\b+\g)z\bar w) \\
\end{array}%
\right)(1-z\bar w)^{-\a-\b-\g-2},$$ for $z,w\in\D.$ We observe that
$$G_1(z,w)=K(0,0)^{1/2}K(z,0)^{-1}K(z,w)K(0,w)^{-1}K(0,0)^{1/2}$$ 
and $G_1(z,0)=I$, as is to be expected.  
The multiplication operator corresponding to $G_1$,  which we denote 
by $M_1$, is \u\ to $M^{(\a,\b+\g)}$ by Lemma \ref{kur}. Hence it is in $\mathrm B_2(\D)$ by
\cite[Proposition 3.6]{DMV}. Since both homogeneity and irreducibility are invariant under
\ue\ it follows, by an easy application of Lemma \ref{kur}, Theorem \ref{tirr}~ and ~Theorem
\ref{homo} that $M_1^*$ is a irreducible homogeneous operator in $\mathrm B_2(\D)$.  
Irreducibility of $M_1^*$ also follows from Proposition
\ref{irr}. Let $M_2$ be the multiplication operator on the Hilbert space of scalar
valued holomorphic functions with reproducing kernel $G_2$. Again, $M_2^*$ is in $\mathrm
B_1(\D).$ The operator $M_2$ is  irreducible by \cite[corollary 1.19]{C-D}.
Homogeneity of $M_2^*$ was first established in \cite{GM}, see also \cite{wil}. An alternative proof is obtained
when we observe that $\Gamma :\mbox{M\"{o}b}\times \D\lo\C$, where
$\Gamma_{\p}(z)=\big((\p)^\prime(z)\big)^{\frac{\a+\b+\g}{2}+1}$ is a cocycle such
that $G_2(z,w)=\Gamma_{\p}(z)G_2\big(\p(z),\p(w)\big)\ov{\Gamma_{\p}(z)}$ for
$z,w\in\D,\v\in$M\"{o}b. Now we conclude that ${M^*}$ is homogeneous as it is
\u\ to the direct sum of two homogeneous operators. Also $M^*$ is in
$\mathrm B_3(\D)$ being the direct sum of two operators from the Cowen-Douglas
class.
\end{proof}

\end{document}